\documentclass[11pt,reqno]{amsart}

\usepackage[margin=1.1in]{geometry}
\usepackage[T1]{fontenc}
\usepackage{lmodern}
\usepackage{amsmath,amssymb,amsthm,mathtools,mathrsfs,tikz}
\usepackage{xcolor}
\usepackage{microtype}
\usepackage[hidelinks]{hyperref}
\usepackage{orcidlink}
\usepackage{tikz-cd}
\definecolor{vegasgold}{rgb}{0.77, 0.7, 0.35}
\definecolor{darkgoldenrod}{rgb}{0.72, 0.53, 0.04}
\definecolor{gold(metallic)}{rgb}{0.83, 0.69, 0.22}

\hypersetup{
 colorlinks=true,
 linkcolor=darkgoldenrod,
 filecolor=brown,
 urlcolor=gold(metallic),
 citecolor=darkgoldenrod,
}

\usepackage{enumitem}
\usepackage[nameinlink,capitalize]{cleveref}

\numberwithin{equation}{section}

\newtheorem{theorem}{Theorem}[section]
\newtheorem{proposition}[theorem]{Proposition}
\newtheorem{lemma}[theorem]{Lemma}

\theoremstyle{definition}
\newtheorem{definition}[theorem]{Definition}

\newcommand{\Z}{\mathbb Z}
\newcommand{\cO}{\mathcal O}

\newcommand{\Fp}{\mathbb{F}_p}

\author[J.~Kalita]{Jayashree Kalita}
\address{Department of Mathematics,
	Vanderbilt University, Nashville, TN 37240}
\email{jayashree.kalita@vanderbilt.edu}

\author[A.~Ray]{Anwesh Ray}
\address{Chennai Mathematical Institute, H1, SIPCOT IT Park, Siruseri, Kelambakkam 603103, Tamil Nadu, India.}
\email{anwesh@cmi.ac.in}

\DeclareMathOperator{\ord}{ord}

\DeclareMathOperator{\Gal}{Gal}

\newcommand{\F}{\mathbb F}

\newcommand{\PP}{\mathbb P}
\newcommand{\Cart}{\mathcal C}

\DeclareMathOperator{\Res}{Res}

\DeclareMathOperator{\Aut}{Aut}

\begin{document}
\raggedbottom
\title[The local lifting conjecture for $D_{p^n}$]
{The local lifting conjecture for dihedral groups of order $2p^n$}
\date{November 22, 2026}

\begin{abstract}
Let $k$ be an algebraically closed field of odd characteristic $p$, and
let $D_{p^n}$ denote the dihedral group of order $2p^n$. A local
$D_{p^n}$-action is a faithful action of $D_{p^n}$ on the formal power
series ring $k[[z]]$ by continuous $k$-automorphisms. The strong Oort
conjecture of Chinburg, Guralnick and Harbater predicts that every such
action lifts to characteristic zero. This is also the dihedral case of
the local lifting conjecture of Obus.

We prove this assertion for $n=2$ and every odd prime $p$. More
precisely, given a local $D_{p^2}$-action on $k[[z]]$, we construct
an action on $R[[Z]]$ by $R$-automorphisms whose reduction modulo the
maximal ideal of $R$ is the given action, where $R$ is a complete
discrete valuation ring finite over the Witt ring $W(k)$.

For arbitrary $n$, let $C_{p^n}$ be the cyclic subgroup of index two
in $D_{p^n}$, and let $u_1,\ldots,u_n$ be the upper ramification
jumps of $k((z))/k((z))^{C_{p^n}}$. We prove that the action lifts
if $u_i\equiv-1\pmod{2p}$ for every $i$. When $p=3$, we also prove
lifting under the hypotheses $u_i\equiv5\pmod6$ for $2\leq i<n$
and $3\nmid u_n$.
\end{abstract}
\hypersetup{pdftitle={Lifting dihedral actions},pdfauthor={Jayashree Kalita and Anwesh Ray}}
\maketitle

\section{Introduction}

\subsection{Motivation}

Let $k$ be an algebraically closed field of characteristic $p>0$, and
let $G$ be a finite group acting faithfully on a smooth projective
connected curve $Y$ over $k$. The lifting problem asks whether $Y$,
together with its $G$-action, is the reduction of a curve with
$G$-action in characteristic zero. More precisely, let $W(k)$
denote the ring of Witt vectors of $k$. We seek a complete discrete
valuation ring $R$ finite over $W(k)$ and a smooth projective
curve $\mathcal Y$ over $R$ with a $G$-action such that
$\mathcal Y\otimes_R k\simeq Y$ as curves with $G$-action
\[
\begin{tikzcd}[row sep=large, column sep=large]
Y \arrow[r, hook, "\iota"] \arrow[d]
  & \mathcal Y \arrow[d] \\
\operatorname{Spec} k \arrow[r, hook]
  & \operatorname{Spec} R .
\end{tikzcd}
\]
This question can be studied one ramification point at a time.
For $y\in Y(k)$, let $G_y$ be its stabilizer. Choosing a local
parameter $z$ at $y$ gives an action of $G_y$ on
$\widehat{\cO}_{Y,y}\simeq k[[z]]$. By the local--global
principle \cite[Theorem~3.1]{ObusSurvey}, the action on $Y$ lifts
over $R$ if and only if each of these local actions lifts over
$R$. It suffices to take one point in each orbit with nontrivial
stabilizer. It is easy to see that there are only finitely many such orbits. Thus in particular, the
local lifts may be taken over a common finite extension of $W(k)$.

Let $k$ be an algebraically closed field of characteristic $p>0$.
A finite group $G$ admits a faithful action on $k[[z]]$ by
$k$-algebra automorphisms if and only if $G$ is isomorphic to a
semidirect product $P\rtimes C_m$, where $P$ is a finite $p$-group
and $C_m$ is cyclic of order $m$ prime to $p$
\cite[Section~2]{CGH}. Such a group $G$ is called a \emph{local Oort
group} for $p$ if, over every algebraically closed field of
characteristic $p$, every local $G$-action lifts to characteristic
zero. The local Oort conjecture, proved by Obus and Wewers \cite{OW}
and Pop \cite{Pop}, asserts that every cyclic group has this property.

Suppose that $p$ is odd. If $G$ is a local Oort group for $p$, then
either $G$ is cyclic or $G$ is a dihedral group of order $2p^n$ for
some $n\geq 1$ \cite[Proposition~1.2]{DDKOT}. Write $D_{p^n}$ for
the dihedral group of order $2p^n$. The conjecture of Obus predicts
that $D_{p^n}$ is a local Oort group for every $n\geq 1$
\cite[Conjecture~1.9 and Remark~1.10]{Obus}.

One necessary condition for a local $G$-action to lift is the
vanishing of the Katz--Gabber--Bertin obstruction. For $G=D_{p^n}$, this condition has a simple form. Let
$u_1$ be the first upper ramification jump of the cyclic extension
$k((z))/k((z))^{C_{p^n}}$. By \cite[Proposition~1.6]{Obus}, the
obstruction vanishes if and only if $u_1$ is odd. Since $u_1$ is
odd for every local $D_{p^n}$-action \cite[Remark~1.7]{DDKOT}, every
such action satisfies this necessary condition. The conjecture
asserts that each of these actions does in fact lift.

Throughout the remainder of the paper, $p$ is odd. We write
\[
 D_{p^n}=\langle \sigma,\tau\mid
 \sigma^{p^n}=\tau^2=1,\ \tau\sigma\tau^{-1}=\sigma^{-1}\rangle.
\]
Given a local $D_{p^n}$-action on $k[[z]]$, set $L=k((z))$ and
$E=L^{\langle\sigma\rangle}$. We denote the upper ramification
jumps of the cyclic extension $L/E$ by $u_1,\ldots,u_n$. All
references to jumps in this paper are to this extension of degree
$p^n$.

Bouw and Wewers proved that every local $D_p$-action lifts
\cite{BW}. For $n=2$, Obus proved that $D_9$ is a local Oort
group and that, for every odd $p$, a local $D_{p^2}$-action
lifts if $u_1\equiv1\pmod p$
\cite[Corollaries~1.18 and~1.19]{Obus}. Dang, Das, Karagiannis,
Obus and Thatte proved that $D_{25}$ and $D_{27}$ are local Oort
groups \cite[Theorem~1.3]{DDKOT}. For arbitrary $n$, Obus proved
lifting under the congruences $u_i\equiv p^{i-1}\pmod{2p}$ for
$1\leq i\leq n$ \cite[Corollary~1.20]{Obus}. Kontogeorgis and Terezakis exhibited
a local $D_{125}$-action, with $p=5$ and lower jumps $9,189,4689$, and claimed
that it does not lift to characteristic zero (see \cite{KT}). Dang and Vasiu show that their construction is incorrect (see
\cite[Theorem~1.2 and Example~14.5]{DV}).

\subsection{Main results}

Our first result treats all odd primes $p$ when $n=2$, without
any restriction on the two jumps.

\begin{theorem}\label{thm-main}
Let $k$ be an algebraically closed field of odd characteristic
$p$. Every local $D_{p^2}$-action on $k[[z]]$ lifts to an action
on $R[[Z]]$ for some complete discrete valuation ring $R$ finite
over $W(k)$.
\end{theorem}

The construction of Bouw and Wewers \cite{BW} is essential to the
proof of Theorem~\ref{thm-main}. In the bounded case, we begin with
their lift of the $D_p$-quotient and construct the remaining cyclic
cover of degree $p$ over it. The local--global principle
therefore gives the global consequence of Theorem~\ref{thm-main}. In greater detail,
every faithful action of $D_{p^2}$ on a smooth projective connected
curve over $k$ lifts to characteristic zero.

For arbitrary $n$, we give two sufficient conditions for lifting.
The first applies in every odd characteristic.

\begin{theorem}\label{thm-arbitrary-length}
Let $k$ be an algebraically closed field of odd characteristic
$p$, and let $n\geq1$. Suppose that a local $D_{p^n}$-action
over $k$ has cyclic upper jumps satisfying
$u_i\equiv-1\pmod{2p}$ for $1\leq i\leq n$. Then the action
lifts over a complete discrete valuation ring finite over $W(k)$.
\end{theorem}

For example, Theorem~\ref{thm-arbitrary-length} applies to every
action with jumps $u_i=2p^i-1$ for $1\leq i\leq n$. In
characteristic three we can allow either possible residue class
of the first jump modulo six, and impose a congruence modulo six
only on the intermediate jumps.

\begin{theorem}\label{thm-characteristic-three}
Let $k$ be an algebraically closed field of characteristic three,
and let $n\geq1$. Suppose that a local $D_{3^n}$-action over $k$
has cyclic upper jumps satisfying
\[
 u_i\equiv5\pmod6\ \text{ for }2\leq i<n
 \quad\text{and}\quad
 3\nmid u_n.
\]
Then the action lifts over a complete discrete valuation ring
finite over $W(k)$.
\end{theorem}

The proof of Theorem~\ref{thm-main} starts with the $D_p$-quotient
of the given $D_{p^2}$-action. Bouw and Wewers constructed a lift of
this quotient to characteristic zero \cite{BW}. We use their lift
as the first step and construct the remaining cover of degree $p$
above it. The main issue is where to place the new branch points.
Their number is determined by the ramification jumps, but the jumps
alone do not determine the original action. We first arrange the
points in pairs exchanged by the reflection so as to obtain the
required ramification. We then vary their positions so that reduction
recovers the given equation, including its coefficients.

For Theorem~\ref{thm-arbitrary-length}, the congruences on the jumps
allow us to move all the wild ramification away from the origin to
points with cyclic stabilizers. The cyclic Oort theorem then suffices.
For Theorem~\ref{thm-characteristic-three}, we instead construct an
explicit polynomial satisfying the lifting criterion of Obus.
\subsection*{Use of artificial intelligence}
The authors used ChatGPT 6.1 Sol as an exploratory tool to study related papers and for verifying references.

\section{Preliminaries}
\subsection{Local actions and ramification}\label{sec-local}

Let $k$ be algebraically closed of odd characteristic $p$. We give a
complete local ring its maximal ideal topology.

\begin{definition}\label{def-local}
Let $G$ be a finite group. A \emph{local $G$-action} on $k[[z]]$ is
an injective homomorphism $\bar\rho$ from $G$ to
$\Aut_k^{\mathrm{cont}}(k[[z]])$. Let $R$ be a complete discrete
valuation ring of mixed characteristic $(0,p)$ with uniformizer $\pi$
and residue field $k$. A \emph{lift} of $\bar\rho$ over $R$ is a
homomorphism $\rho$ making the following diagram commute.
\[
\begin{tikzcd}[row sep=large,column sep=large]
 &\Aut_R^{\mathrm{cont}}(R[[Z]])\arrow[d,"\mathrm{red}"]\\
 G\arrow[ur,dashed,"\rho"]\arrow[r,"\bar\rho"']&\Aut_k^{\mathrm{cont}}(k[[z]]) .
\end{tikzcd}
\]
Here reduction sends an automorphism to its action modulo $\pi$, with
the image of $Z$ identified with $z$.
\end{definition}

The local Oort conjecture asserts that every local cyclic action lifts
over a complete discrete valuation ring finite over $W(k)$
\cite[Conjecture~1.3]{OW}.

\begin{theorem}[Local--global principle {\cite[Theorem~3.1]{ObusSurvey}}]\label{thm-local-global}
Let $G$ act faithfully on a smooth projective connected curve $Y$ over
$k$, and let $R$ be a complete discrete valuation ring finite over
$W(k)$. The action lifts to a smooth projective curve over $R$ if and
only if, for each $y\in Y(k)$, the action of its stabilizer $G_y$ on
$\widehat{\cO}_{Y,y}$ lifts over $R$. It suffices to take one point in
each orbit with nontrivial stabilizer.
\end{theorem}

If $G$ is cyclic, each $G_y$ is cyclic. There are only finitely many
orbits to consider, so the local lifts can all be taken over a common
finite extension of $W(k)$. Theorem~\ref{thm-local-global} then gives a
global lift. Conversely, the Katz--Gabber theorem realizes a given
local cyclic action as the completion of a cover of the projective line
\cite[Theorem~1.4.1 and Corollary~1.4.12]{Katz}. Completing a lift of
this cover gives a lift of the local action.

Let $L=k((z))$ and $E_0=L^G$. The extension $L/E_0$ attached to a local
$G$-action is totally ramified of degree $|G|$, since $k$ is
algebraically closed. For $i\geq0$ let $G_i$ denote its $i$th lower
ramification group, and write $d(L/E_0)$ for its different exponent.
With $\ord_z(0)=+\infty$, Hilbert's
different formula \cite[Chapter~IV, \S1, Proposition~4]{Serre} gives
\begin{equation}\label{pre-hilbert}
d(L/E_0)=\sum_{i\geq0}(|G_i|-1).
\end{equation}
If $G_i\ne G_{i+1}$,
then $i$ is a lower jump and $\varphi(i)$ is the corresponding upper
jump, where \[\varphi(r)=\int_0^r |G_{\lceil t\rceil}|/|G_0|\,dt.\]
Upper numbering is compatible with passage to quotients
\cite[Chapter~IV, \S3]{Serre}.

Suppose $L/E_0$ is cyclic of degree $p^2$ with upper jumps $a<M$.
The lower ramification groups have order $p^2$ for $0\le i\le a$
and order $p$ for $a<i\le a+p(M-a)$, since the slope of $\varphi$
between the two jumps is $1/p$. Denote the different exponent by $d(a,M)$.
By \eqref{pre-hilbert}, it equals $(a+1)(p^2-1)+p(M-a)(p-1)$.
In terms of the two upper jumps this gives
\begin{equation}\label{eq-different}
 d(a,M)=(p-1)(a+1)+(p^2-p)(M+1).
\end{equation}
To compare with a cover in characteristic zero, write
\begin{equation}\label{eq-budget}
 d(a,M)=(a+1)(p^2-1)+(M-a)(p^2-p).
\end{equation}
The right hand side is the different degree of a cyclic cover of
degree $p^2$ in characteristic zero with $a+1$ branch points of index
$p^2$ and $M-a$ of index $p$. Indeed, a branch point of index $p$ has
$p$ points above it, each with different exponent $p-1$.

\subsection{Witt equations for a dihedral action}\label{pre-asw}

Set $E=k((x))$ and $G_E=\Gal(E^{\mathrm{sep}}/E)$. For $f\in E$, a
root $y$ of $y^p-y=f$ determines a character $g\mapsto g(y)-y$ with
values in $\F_p$. This character depends only on the class of $f$
modulo elements of the form $b^p-b$, since $y+b$ is a root of the
equation with right hand side $f+b^p-b$ for $b\in E$.

For a cyclic extension of degree $p^n$, set $\wp=F-1$ on the group
$W_n(E)$ of Witt vectors of length $n$. Here $F$ raises each
coordinate to the $p$th power, and subtraction is Witt subtraction.
The Artin--Schreier--Witt sequence gives the isomorphism
\begin{equation}\label{eq-witt}
 W_n(E)/\wp W_n(E)\ \xrightarrow{\ \sim\ }\
 \operatorname{Hom}_{\mathrm{cont}}(G_E,\Z/p^n\Z),
\end{equation}
given by $[f]\longmapsto\bigl(g\longmapsto g(y)-y\bigr)$, where $y$ is any solution of $\wp(y)=f$ \cite[Section~2]{Obus}.
Indeed $g(y)-y\in W_n(\F_p)=\Z/p^n\Z$. The fixed field of the kernel
is $E(y_1,\ldots,y_n)$, and the degree-$p$ quotient of a character of
order $p^n$ is given by $y_1^p-y_1=f_1$. We call the character in
\eqref{eq-witt} the \emph{Witt character} of the cyclic extension once
a generator of its Galois group has been chosen.

For a local $D_{p^n}$-action, set $L=k((z))$ and
$E=L^{\langle\sigma\rangle}$. Let $t$ be a uniformizer of $E$. The coefficient of $t$
in $\tau(t)$ is $1$ or $-1$, since $\tau^2=1$. Suppose it is $1$.
As $\tau$ acts nontrivially on $E$, we may write
\[\tau(t)=t+ct^m+O(t^{m+1})\] with $m\ge2$ and $c\ne0$. Then
\[\tau^2(t)\equiv t+2ct^m\pmod{t^{m+1}},\] a contradiction.
Set $x=(t-\tau(t))/2$. Then $x\equiv t\pmod{t^2}$ and $\tau(x)=-x$.
Setting $s=x^2$, we obtain the tower
\[
 L^{D_{p^n}}=k((s))\ \subset\ E=k((x))\ \subset\ L
\]
of successive degrees $2$ and $p^n$. Call a polynomial in $x^{-1}$
\emph{standard} if its nonzero terms have positive degree prime to $p$.
This includes the zero polynomial. A function $h$ is \emph{odd} if
$h(-x)=-h(x)$.

\begin{lemma}\label{arb-standard-witt}
Every local $D_{p^n}$-action has an equation $\wp(y)=f$ in which
$f_1,\ldots,f_n$ are odd standard polynomials in $x^{-1}$ and
$\sigma(y)=y+1$, $\tau(y)=-y$. Set $d_i=\deg_{x^{-1}}f_i$, with
$\deg 0=-\infty$. Then the cyclic upper jumps satisfy
\begin{equation}\label{arb-jump-formula}
 u_1=d_1\quad\text{and}\quad u_i=\max\{pu_{i-1},d_i\}\quad \text{for}\quad 2\leq i\leq n.
\end{equation}
\end{lemma}
\begin{proof}
Choose the character sending $\sigma$ to $1$ and let $f$ represent its
class in \eqref{eq-witt}. Since $\tau\sigma\tau^{-1}=\sigma^{-1}$,
we have $[\tau f]=-[f]$. Multiplication by $2$ is invertible on
$W_n(E)$, so $(f-\tau f)/2$ represents the same class. Replacing $f$
by this vector gives $\tau f=-f$. For odd $p$,
negating each coordinate negates each Witt polynomial, since all
powers occurring in these polynomials are odd. Thus Witt negation is
coordinatewise and $f_i(-x)=-f_i(x)$ for every $i$.

We make $f_i$ standard in order of increasing $i$. Subtract
$\wp(0,\ldots,0,b,0,\ldots,0)$ with $b$ in the $i$th coordinate.
This changes $f_i$ by $-(b^p-b)$ and fixes the preceding coordinates.
For the regular part $h\in xk[[x]]$, take
$b=-\sum_{j\geq0}h^{p^j}$, so that $b^p-b=h$. For a pole term
$cx^{-pm}$, take $b=c^{1/p}x^{-m}$. This replaces the term by one of
smaller pole order. Starting with the largest pole order divisible by
$p$ gives a standard polynomial after finitely many such subtractions.
In both cases $b$ is odd. Thus the new Witt vector still satisfies
$\tau f=-f$, including in the coordinates yet to be treated.

Choose a solution $y$ of $\wp(y)=f$ in $W_n(L)$. Set
$c=\tau(y)+y$. Since $\wp(c)=\tau(f)+f=0$, we have
$c\in W_n(\F_p)$, so $\tau$ and $\sigma$ fix $c$. The solution
$y'=y-c/2$ then satisfies $\tau(y')=-y'$ and $\sigma(y')=y'+1$.
Formula~\eqref{arb-jump-formula} is the upper jump formula for a
standard Witt vector \cite[Section~2, equation~(2.1)]{Obus}. The
first coordinate is nonzero since the character has order $p^n$.
\end{proof}

Since each $f_i$ is a polynomial in $x^{-1}$, we may regard
$\wp(y)=f$ as an equation over $k(x)$. Choose the solution
$y=(y_1,\ldots,y_n)$ in the given local field $L$, and set
$F=k(x)(y_1,\ldots,y_n)$. Let $Y$ be the smooth projective curve
with function field $F$. The inclusions $k(s)\subset k(x)\subset F$
give maps
\[
 Y\longrightarrow\PP^1_x\longrightarrow\PP^1_s,
\]
where $s=x^2$. Consider first the map $Y\to\PP^1_x$. Completing $F$ at the point
determined by our choice of $y$ gives
$k((x))(y_1,\ldots,y_n)=L$. Since $[L:k((x))]=p^n$ and the Witt
equation has degree $p^n$, we have $[F:k(x)]=p^n$. Thus there is
a unique point $P$ above $x=0$, and its ramification index is $p^n$.
At every other point of the $x$-line, including infinity, each
$f_i$ is regular. The successive Artin--Schreier equations are
monic with derivative $-1$, so the cover is unramified there.

Now consider the composite map $Y\to\PP^1_s$. Translation by the
Witt vector $1$ gives an automorphism $\sigma$ fixing $x$. Since
$f(-x)=-f(x)$, the substitution $\tau(x,y)=(-x,-y)$ also preserves
the equation. Both automorphisms fix $s=x^2$, and
$\tau\sigma\tau^{-1}=\sigma^{-1}$. They therefore generate the
dihedral group of order $2p^n$, which is the degree of the
composite map.

The map $x\mapsto x^2$ has ramification index two at zero and at
infinity and is unramified elsewhere. Hence the composite has
index $2p^n$ at $P$ and index two at each point above $s=\infty$.
It has no other branch points. Its completed extension at $P$ is
the given dihedral extension $L/k((s))$.

\begin{lemma}\label{pre-normalized-witt}
For $n=2$, one may choose $x$ with $\tau(x)=-x$ so that the Witt
representative is $(x^{-a},f_2)$ with $f_2$ odd and standard.
The cyclic upper jumps are $a$ and $M=\max\{pa,\deg_{x^{-1}}f_2\}$.
Both are odd and $p\nmid a$.
\end{lemma}
\begin{proof}
In the odd parameter $t$ of Lemma~\ref{arb-standard-witt}, write
$f_1=ct^{-a}U(t)$ with $U(0)=1$. Since $f_1$ is odd, $a$ is odd
and $U(-t)=U(t)$. Choose $b\in k$ with $b^a=c$. As $p\nmid a$,
Hensel's lemma gives a unique $V\in1+tk[[t]]$ with $V^a=U$.
Uniqueness implies $V(-t)=V(t)$, so $x=t/(bV)$ is odd and
$f_1=x^{-a}$. Making the second coordinate standard as in
Lemma~\ref{arb-standard-witt} leaves the first coordinate unchanged.
Formula~\eqref{arb-jump-formula} now gives the assertion.
\end{proof}

The degree-$p$ quotient thus has equation $w^p-w=x^{-a}$ with
reflection $(x,w)\mapsto(-x,-w)$. We use the following theorem to
lift the resulting dihedral action.

\begin{theorem}[Bouw--Wewers {\cite[Theorem~1.2]{BW}}]\label{ht-classical-theorem}
Every local $D_p$-action over $k$ lifts to characteristic zero.
\end{theorem}

For the proof in degree $p^2$, we shall use the lift given by
\cite[Lemma~4.2]{Obus}. We shall take a note of its branch points and their positions
in Section~\ref{sec-analytic}.

\section{The Cartier equation and its coefficients}\label{sec-cartier}

Let $a<p$ be the first cyclic upper jump of the given $D_{p^2}$-action.
Recall that $a$ is odd. We have chosen $x$ so that the reflection
sends $x$ to $-x$ and the Witt vector is $(x^{-a},f_2(x))$. Let
$Y_2$ be the cover of the $x$-line defined by this vector. Its
quotient $Y_1$ of degree $p$ is given by $w^p-w=x^{-a}$. We fix a
lift $\mathcal Y_1$ of $Y_1$, together with its reflection, as in
\cite[Theorem~1.2]{BW} and \cite[Lemma~4.2]{Obus}.

Let $R$ be the coefficient ring of this lift, with maximal ideal
$\mathfrak m$, residue field $k$ and fraction field $E$. The diagram
below shows what remains to be constructed. Each term on the left
is obtained from the term on its right by reduction modulo
$\mathfrak m$.
\[
\begin{tikzcd}[row sep=large,column sep=large]
Y_2 \arrow[d,"p"'] \arrow[r,dashed] &
\mathcal Y_2 \arrow[d,dashed,"p"] \\
Y_1 \arrow[d,"p"'] \arrow[r,hook] &
\mathcal Y_1 \arrow[d,"p"] \\
\mathbb P^1_k \arrow[r,hook] & \mathbb P^1_R .
\end{tikzcd}
\]
Here $p$ denotes the degree of each vertical map. We must choose
the new branch points of $\mathcal Y_2$ so that its reduction gives
the prescribed second Witt coordinate $f_2$.

Let $T$ be the variable on the base $\mathbb P^1_R$, chosen so that
$T$ reduces to $x$ and the reflection sends $T$ to $-T$. Normalize
the valuation on $E$ by $v(p)=1$ and set $r=1/((p-1)a)$. After a
finite extension of $E$, choose $\lambda$ with $v(\lambda)=r$.
To distinguish branch points $T=b$ with $v(b)=r$, we reduce
$U=T/\lambda$ instead of $T$. Write $t$ for the reduction of $U$.

We now describe the differential contributed by the fixed cover
$\mathcal Y_1$. After adjoining the $p^2$-th roots of unity to $E$,
write its equation as $Y^p=F(T)$. By
\cite[Proposition~5.17 and Lemma~6.1]{OW}, we may divide $Y$ by a
rational function of $U$ so that the equation becomes
$Y^p=1+pA(U)$, where $A$ has a reduction $h\in k(t)$. Explicitly,
if $A=P/Q$ with $P,Q\in R[U]$ and $\overline Q\ne0$, then
$h=\overline P/\overline Q$. The differential we use is the
derivative of this function. For the chosen lift it is
\[
 \omega_1=dh=h'(t)\,dt=c\,\frac{dt}{t^{a+1}}
\]
for some $c\in k^\times$. Thus $c$ is fixed before we choose the
new branch points.

Write the proposed cover of degree $p^2$ as $Z^{p^2}=B(T)$. A branch
point is new if the degree-$p$ quotient is unramified there.
Among the new branch points, write those with $v(b)=r$ as
$\pm b_1,\ldots,\pm b_n$. Set $x_i=\overline{b_i/\lambda}$.
At $b_i$ the integer $m_i=\operatorname{ord}_{T=b_i}B$ is divisible
by $p$ but not by $p^2$. Let $e_i\in\mathbb F_p^\times$ be the
reduction of $m_i/p$ modulo $p$. Since the reflection acts by
inversion, the exponent at $-b_i$ gives $-e_i$.

We can now state the condition on the unknown positions $x_i$.
If the cover has good reduction, the $2n$ elements $\pm x_i$ are
distinct. Moreover, by \cite[Proposition~7.9 and Corollary~7.11]{Obus},
the even polynomial
\[
 f(t)=-\frac{1}{c}\prod_{i=1}^n
              \left(1-\frac{t^2}{x_i^2}\right)
\]
must satisfy
\[
 \frac{dt}{t^{a+1}f(t)}
 =
 \sum_{i=1}^n e_i
 \left(\frac{dt}{t-x_i}-\frac{dt}{t+x_i}\right)-\omega_1 .
\]
Thus the roots of $f$ give the reductions of the new branch points
at this valuation. The displayed identity is the equation that
these roots must satisfy.

Set $\omega_f=dt/(t^{a+1}f(t))$. The Cartier operator $\Cart$ fixes
each $dt/(t-b)$ and sends each derivative $dh$ to zero. Since
$e_i\in\mathbb F_p$ and $\omega_1=dh$, the identity above gives
\[
 \Cart(\omega_f)-\omega_f=c\,\frac{dt}{t^{a+1}}.
\]
For fixed exponents $e_i$ and the degrees considered below, we shall
prove that only finitely many tuples $(x_1,\ldots,x_n)$ occur.
We shall also give an explicit $f$ in each of the two degrees
needed for the lift. We begin by defining $\Cart$ on a rational
differential $H(t)\,dt$ with $H\in k(t)$.
\begin{lemma}\label{lem-rational-frobenius}
Let $k$ be a perfect field of characteristic $p>0$. Every $H\in k(t)$
can be written uniquely as $H=\sum_{i=0}^{p-1}g_i^pt^i$, where
$g_i\in k(t)$.
\end{lemma}

\begin{proof}
Suppose first that $H$ is a polynomial. Group its terms according
to their exponents modulo $p$. A term $ct^{pm+i}$ can be written as
$(c^{1/p}t^m)^pt^i$, since $k$ is perfect. Adding the terms in
each group gives the required expression. For a rational function
$H=A/B$, apply this construction to the polynomial $AB^{p-1}$ and
divide by $B^p$.

For uniqueness, suppose $\sum_i h_i^pt^i=0$ and choose a common
denominator $D$ for the $h_i$. Every term of $(Dh_i)^pt^i$ has
degree congruent to $i$ modulo $p$. Thus no terms from distinct
summands can cancel in $\sum_i(Dh_i)^pt^i=0$. It follows that
$h_i=0$ for all $i$. Applying this to the difference of two
expressions for $H$ proves uniqueness.
\end{proof}

With $H$ as in Lemma~\ref{lem-rational-frobenius}, we define
$\Cart(H\,dt)=g_{p-1}\,dt$. It follows from the definition that
$\Cart$ is additive and that $\Cart(A^p\omega)=A\Cart(\omega)$ for
$A\in k(t)$. In particular, for $b\in k$ and $j\in\Z$, we have
\begin{equation}\label{cart-eq-monomial}
 \Cart(b^pt^j\,dt)=
 \begin{cases}
 b\,t^{(j+1)/p-1}\,dt & \text{if }p\mid j+1,\\
 0 & \text{otherwise.}
 \end{cases}
\end{equation}

For the remainder of the section, $k$ is algebraically closed of
odd characteristic $p$. A differential is called \emph{exact} if
it equals $dA$ for some $A\in k(t)$, and \emph{logarithmic} if it
equals $dA/A$ for some $A\in k(t)^\times$. At a finite point $b$,
we expand the differential in the variable $u=t-b$. Its residue at
$b$ is the coefficient of $u^{-1}du$. At infinity we use $u=1/t$.

\begin{lemma}\label{cart-lem-exact-log}
A rational differential $\omega$ is exact if and only if
$\Cart(\omega)=0$. The following are equivalent.
\begin{enumerate}
\item $\omega$ is logarithmic.
\item Every pole of $\omega$ is simple and its residue belongs to $\Fp$.
\item $\Cart(\omega)=\omega$.
\end{enumerate}
\end{lemma}

\begin{proof}
Write $\omega=(\sum_{i=0}^{p-1}g_i^pt^i)\,dt$. If $g_{p-1}=0$,
then $\omega$ is the derivative of
$\sum_{i=0}^{p-2}g_i^pt^{i+1}/(i+1)$, since $d(g_i^p)=0$.
Conversely, if $A=\sum_{i=0}^{p-1}h_i^pt^i$, then
$dA=\sum_{i=1}^{p-1}i h_i^pt^{i-1}\,dt$. There is no term of
index $p-1$ in this expression, so $\Cart(dA)=0$.

Write $A=b_0\prod_b(t-b)^{m_b}$, with $m_b\in\Z$. Then
$dA/A=\sum_b m_b\,dt/(t-b)$. Its poles are simple, and its
residues at the finite points and at infinity are the images in
$\Fp$ of the integers $m_b$ and $-\sum_bm_b$, respectively. This
proves that (1) implies (2). Conversely, suppose (2) holds.
For each finite pole $b$, choose an integer $m_b$ representing its
residue. Subtracting $\sum_bm_b\,dt/(t-b)$ from $\omega$ leaves a
differential $P(t)\,dt$ with $P$ a polynomial and at most a simple
pole at infinity. If $P\ne0$, substitution of $t=1/u$ gives a
pole of order $\deg P+2$. Hence $P=0$, proving (1).

To prove the equivalence with (3), we may use
\eqref{cart-eq-monomial} in the variables $u=t-b$ and $u=1/t$.
Indeed, $u^i\,du$ is exact for $0\le i<p-1$. For $u=t-b$,
expansion of $(t-b)^{p-1}dt$ gives $\Cart(u^{p-1}du)=du$.
For $u=1/t$, the same identity follows from
$u^{p-1}du=-t^{-p-1}dt$. The rule $\Cart(A^p\omega)=A\Cart(\omega)$
now gives \eqref{cart-eq-monomial} for every integral power of $u$.
In particular, $\Cart(du/u)=du/u$. Applying this to the expression
for $dA/A$ proves that (1) implies (3).

Finally, suppose $\Cart(\omega)=\omega$. At a pole of order
$m\ge2$, formula~\eqref{cart-eq-monomial} gives a pole of order at
most $1+\lfloor(m-1)/p\rfloor<m$ after applying $\Cart$, if a pole
remains at all. Thus every pole of $\omega$ is simple. At a simple
pole of residue $b$, the residue after applying $\Cart$ is
$b^{1/p}$. Hence $b=b^{1/p}$, so $b\in\Fp$. This proves (2).
\end{proof}

Fix $c\in k^\times$. We consider nonzero even polynomials $f\in k[t]$
for which
\begin{equation}\label{eq-cartier}
 \Cart(\omega_f)-\omega_f=c.
\end{equation}
Here $a$ and $c$ are fixed, and the coefficients of $f$ are the
unknowns. We consider the degrees $N=(p-1)a-2d\ge0$, where
$0\le d<p$. The two cases used in the lifting argument are $d=0$
and $d=a$.

\begin{definition}\label{def-ordinary}
An \emph{ordinary datum} for $(a,d)$ is an even polynomial of degree
$N=(p-1)a-2d$ satisfying \eqref{eq-cartier} with $c=a$.
\end{definition}

Let $f$ satisfy \eqref{eq-cartier}, and set
\[
 \eta=\omega_f+c\,\frac{dt}{t^{a+1}}
     =\left(\frac1{f(t)}+c\right)\frac{dt}{t^{a+1}}.
\]
Since $c\,dt/t^{a+1}=d(-(c/a)t^{-a})$ is exact, we have
$\Cart(\eta)=\eta$. By Lemma~\ref{cart-lem-exact-log}, every pole
of $\eta$ is simple and its residue belongs to $\Fp$.

We first consider $t=0$. If $f(0)=0$, then $\eta$ has a pole of
order greater than $a+1$ at zero, which is impossible. Thus
$f(0)\ne0$. The coefficient of $t^{-a-1}dt$ in $\eta$ is
$1/f(0)+c$. Since $a+1>1$, this coefficient must vanish, giving
$f(0)=-1/c$. In fact $\eta$ is regular at zero. Indeed, the
coefficient of $dt$ in $\eta$ is even, so its expansion contains
no term in $t^{-1}$.

Now let $x$ be a root of $f$. Since $x\ne0$, the differential
$c\,dt/t^{a+1}$ is regular at $x$. Thus the multiplicity of $x$
in $f$ is the pole order of $\eta$ at $x$, which is one. Since
$f$ is even, its $N$ roots form $n=N/2$ distinct pairs
$\{x_i,-x_i\}$. Choose one root $x_i$ from each pair. At this root
the residue is
\[
 r_i=\Res_{t=x_i}\eta=\Res_{t=x_i}\omega_f
     =\frac{1}{x_i^{a+1}f'(x_i)}\in\Fp^\times.
\]
The residue at $-x_i$ is $-r_i$, since $f'$ is odd and $a+1$
is even.

We next express \eqref{eq-cartier} in terms of the roots $x_i$
and their residues $r_i$. Set $K=a+N$. The denominator of
$\omega_f$ has degree $K+1$, so its coefficient, expanded in
$t^{-1}$, begins with a nonzero multiple of $t^{-K-1}$.
The differential $\eta$ is regular at zero and at infinity, and
its only poles are the simple poles at $\pm x_i$. Hence
\begin{equation}\label{eq-log}
 \eta=
 \sum_{i=1}^{n}r_i
 \left(\frac{dt}{t-x_i}-\frac{dt}{t+x_i}\right)
 =
 2\sum_{\substack{q\ge1\\q\ {\rm odd}}}
 \left(\sum_{i=1}^n r_i x_i^q\right)t^{-q-1}dt.
\end{equation}
For the first equality, subtract the indicated partial fractions
from $\eta$. The difference is regular everywhere on the projective
line and hence zero. The second equality follows by expanding each
fraction in powers of $t^{-1}$. Comparison with
$\eta=\omega_f+c t^{-a-1}dt$ gives the equations in the following
lemma.

\begin{lemma}\label{lem-finitecartier}
Fix $a,d,c$ as above. The ring obtained by imposing
\eqref{eq-cartier} on the coefficients of an even polynomial $f$
of degree $N$ with $f(0)=-1/c$ is finite dimensional over $k$.
In particular, there are only finitely many such polynomials.

Fix $r_1,\ldots,r_n\in\Fp^\times$. Consider the ordered tuples
$(x_1,\ldots,x_n)$ with $x_i\ne0$ and $x_i\ne\pm x_j$ for
$i\ne j$ satisfying
\begin{equation}\label{eq-moments}
 2\sum_{i=1}^n r_ix_i^q=
 \begin{cases}
 c & \text{if }q=a,\\
 0 & \text{if }q\ne a,
 \end{cases}
 \quad 1\le q<K,\quad q\text{ odd},\quad p\nmid q.
\end{equation}
Their coordinate ring is also finite dimensional over $k$.
Every such tuple gives a polynomial satisfying \eqref{eq-cartier}
with roots $\pm x_i$ and residues $r_i$ at $x_i$.
\end{lemma}

\begin{proof}
Write $f(t)=-1/c+\sum_{j=1}^n b_jt^{2j}$ and set $D=t^{a+1}f(t)$.
Multiplying \eqref{eq-cartier} by $D$ and using
$\Cart(D^p\omega)=D\Cart(\omega)$ gives
\begin{equation}\label{eq-cartierdef}
 \Cart(D^{p-1}dt)=(1+cf(t))dt.
\end{equation}
For $1\le j\le n$, comparison of the coefficients of $t^{2j}dt$
using \eqref{cart-eq-monomial} gives
\begin{equation}\label{eq-frobenius}
 (c b_j)^p=[t^{p(2j+1)-1}]D^{p-1},
\end{equation}
where $[t^m]$ denotes the coefficient of $t^m$. These are all
the coefficient equations. Indeed, both constant terms in
\eqref{eq-cartierdef} vanish, and only even powers of $t$ occur.
To obtain a term of degree $2n+2$ after applying $\Cart$, one
would need a term of degree $p(2n+3)-1$ in $D^{p-1}$. This exceeds
$\deg D^{p-1}=(p-1)(a+1+2n)$ by $2(p-d)>0$.

The right hand side of \eqref{eq-frobenius} has total degree at most
$p-1$ in $b_1,\ldots,b_n$. Thus replacing a factor $b_j^p$ using
this equation lowers the total degree. It follows that the quotient
of $k[b_1,\ldots,b_n]$ by these equations is spanned by the
monomials $\prod_j b_j^{e_j}$ with $0\le e_j<p$. When $n>0$ we
invert $b_n$ to require $\deg f=2n$. The resulting ring is still
finite dimensional over $k$. This proves the first assertion.

For a polynomial satisfying \eqref{eq-cartier}, comparison of the
coefficients of $t^{-q-1}dt$ in \eqref{eq-log} gives
\eqref{eq-moments}. The omitted odd multiples of $p$ give no new
conditions. Indeed, if $pm<K=pa-2d$ with $m$ odd, then $m<a<p$,
and the equation for $m$ gives
$\sum_i r_i x_i^{pm}=(\sum_i r_i x_i^m)^p=0$. There are exactly
$n$ equations in \eqref{eq-moments}. Of the $(K-1)/2$ positive odd
integers below $K$, precisely $(a-1)/2$ are multiples of $p$,
since $p(a-2)<K\le pa$. Their difference is $(K-a)/2=n$.

Conversely, suppose that $(x_1,\ldots,x_n)$ satisfies
\eqref{eq-moments} and the stated inequalities. The only polynomial
of degree $N$ with roots $\pm x_i$ and constant term $-1/c$ is
\[
 f(t)=-\frac1c\prod_{i=1}^n
       \left(1-\frac{t^2}{x_i^2}\right).
\]
Let $R(t)\,dt$ denote the sum of partial fractions in \eqref{eq-log}.
By \eqref{eq-moments} and the preceding observation on multiples
of $p$, we have $R(t)-ct^{-a-1}=O(t^{-K-1})$ at infinity. Set
$P(t)=t^{a+1}f(t)(R(t)-ct^{-a-1})$. The factor $t^{a+1}f(t)$
clears all denominators, so $P$ is a polynomial. Since this factor
has degree $K+1$, the estimate at infinity gives $\deg P\le0$.
Moreover, $P(0)=-cf(0)=1$. Hence $P=1$, or equivalently,
\[
 R(t)\,dt=\omega_f+c\,\frac{dt}{t^{a+1}}.
\]
Cartier fixes the left hand side by Lemma~\ref{cart-lem-exact-log}
and kills $c\,dt/t^{a+1}$. Applying Cartier proves
\eqref{eq-cartier}. Comparing residues at $x_i$ gives the
prescribed value $r_i$.

There are only finitely many possible $f$, and each admits only
finitely many choices of signs and orderings of its roots. Hence
there are only finitely many tuples satisfying \eqref{eq-moments}
and the inequalities in the statement. Their coordinate ring is
obtained from $k[x_1,\ldots,x_n]$ by imposing \eqref{eq-moments}
and inverting $\prod_i x_i\prod_{i<j}(x_i^2-x_j^2)$. By the
Nullstellensatz, its quotient by the ideal of nilpotent elements
is a finite product of copies of $k$, or is zero if there are no
solutions. The ideal of nilpotent elements is finitely generated
and nilpotent. Each successive quotient of its powers is a finite
module over the reduced quotient. Thus the coordinate ring itself
is finite dimensional over $k$. In particular, its completed local
ring at each solution is Artinian.
\end{proof}

\begin{lemma}\label{lem-data}
For $d=0$ and $d=a$, respectively, ordinary data are given by
\begin{equation}\label{eq-explicitdata}
 f_a^{(0)}(t)=\frac{t^{(p-1)a}-1}{a},
 \quad\text{and}\quad
 f_a^{(1)}(t)=-\frac{t^{(p-1)a}-1}{a(t^{2a}-1)}.
\end{equation}
Both polynomials have only simple roots.
\end{lemma}

\begin{proof}
Since $(t^{(p-1)a}-1)/(t^{2a}-1)=\sum_{j=0}^{(p-3)/2}t^{2aj}$,
the two polynomials are even, of respective degrees $(p-1)a$ and
$(p-3)a$, and both have constant term $-1/a$. The polynomial
$t^{(p-1)a}-1$ has simple roots since $p\nmid(p-1)a$. The roots
of the second polynomial are obtained by removing those with
$t^{2a}=1$, so they too are simple.

For either choice of $f$, set $\eta=\omega_f+a\,dt/t^{a+1}$.
At zero, $1/f=-a+O(t^{2a})$, so $\eta=O(t^{a-1})dt$ is regular.
At infinity both terms of $\eta$ are regular. Indeed, in the
coordinate $u=1/t$ the term $t^{-a-1}dt$ becomes $-u^{a-1}du$,
and $\omega_f$ has order $a+\deg f-1\ge0$. It remains to compute
the residues at the roots of $f$.

Let $x$ be such a root and set $z=x^a$. Since $z^{p-1}=1$, we have
$z\in\Fp^\times$. For the first polynomial,
$(f_a^{(0)})'(x)=(p-1)x^{(p-1)a-1}=-1/x$. At a root of the
second polynomial we also have $z^2\ne1$. Differentiating
$f_a^{(1)}=-f_a^{(0)}/(t^{2a}-1)$ at this root gives
$(f_a^{(1)})'(x)=1/(x(x^{2a}-1))$. Thus the residue
$1/(x^{a+1}f'(x))$ is $-z^{-1}$ for the first polynomial and
$z-z^{-1}$ for the second. Both lie in $\Fp^\times$.

It follows from Lemma~\ref{cart-lem-exact-log} that
$\Cart(\eta)=\eta$. Since $\Cart(t^{-a-1}dt)=0$, this proves
\eqref{eq-cartier} with $c=a$. When $p=3$ and $a=1$, the second
polynomial is the constant $-1$ and $\eta=0$, so the same
conclusion holds.
\end{proof}

Each polynomial in \eqref{eq-explicitdata} also gives a solution
for any $c\in k^\times$. Choose $\lambda\in k^\times$ with
$\lambda^a=c/a$ and set $g(t)=\lambda^{-a}f(t/\lambda)$. Then
$g(0)=-1/c$, and its roots are the roots of $f$ multiplied by
$\lambda$. Under the substitution $u=t/\lambda$, the differential
$du/(u^{a+1}f(u))$ becomes $\omega_g$, while $a\,du/u^{a+1}$
becomes $c\,dt/t^{a+1}$. Since Cartier commutes with this
substitution, $g$ satisfies \eqref{eq-cartier} with coefficient $c$.

We now explain why the finiteness in Lemma~\ref{lem-finitecartier}
is useful. The roots $x_i$ belong to $k$. In characteristic zero
we first choose elements reducing to these roots. We then need to
change the chosen elements without changing their reductions, so
that they satisfy the required equations. The following lemma
allows such changes even when the matrix of first derivatives is
singular.

\begin{lemma}\label{lem-formal}
Let $R$ be a complete discrete valuation ring with uniformizer $\pi$,
and set $B=R[[z_1,\ldots,z_m]]$. Let $\phi_1,\ldots,\phi_m\in B$
have zero constant term. Suppose that
\begin{equation}\label{cart-eq-fibre}
 B/(\pi,\phi_1,\ldots,\phi_m)
\end{equation}
has finite positive dimension over $R/\pi R$. For any
$b_1,\ldots,b_m\in\pi R$, there exist a complete discrete valuation
ring $R'$ finite over $R$ and elements $z'_1,\ldots,z'_m$ of its
maximal ideal such that $\phi_i(z'_1,\ldots,z'_m)=b_i$ for all $i$.
\end{lemma}

\begin{proof}
Set $A=R[[y_1,\ldots,y_m]]$ and map $A$ to $B$ by
$y_i\mapsto\phi_i$. Let $I=(\pi,\phi_1,\ldots,\phi_m)B$. Since
$B/I$ is a finite dimensional local algebra, its maximal ideal is
nilpotent. Thus $\mathfrak m_B^h\subseteq I\subseteq\mathfrak m_B$
for some $h$, so $B$ is complete for the powers of $I$.

Choose $e_1,\ldots,e_s\in B$ lifting a basis of $B/I$ over
$R/\pi R$. We show that these elements generate $B$ over $A$.
They generate modulo $I$ by construction. Since
$I^j=\mathfrak m_A^jB$, an error in $I^j$ can be removed modulo
$I^{j+1}$ by adding a linear combination of the $e_i$ with
coefficients in $\mathfrak m_A^j$. The successive coefficients
converge in $A$, and the errors tend to zero in $B$. Hence $B$ is
finite over $A$. Both rings are regular local rings of dimension
$m+1$, and $B/I$ has dimension zero. The flatness criterion
\cite[Tag~00R4]{Stacks} gives that $B$ is flat over $A$. It is
therefore free over $A$ of rank $s=\dim_{R/\pi R}(B/I)>0$.

Now evaluate $y_i$ at $b_i$. The resulting ring $C=B\otimes_A R$
is free of rank $s$ over $R$. Thus $C[1/\pi]$ is a nonzero finite
dimensional algebra over $\operatorname{Frac}(R)$. Its quotient by
a maximal ideal is a finite field extension $E$ of
$\operatorname{Frac}(R)$. The valuation ring $R'$ of $E$ is a
complete discrete valuation ring finite over $R$. The image of $C$
lies in $R'$ since its elements are integral over $R$.

Finally, $C/\pi C=B/I$ because every $b_i$ belongs to $\pi R$.
The image of each $z_i$ in this quotient is nilpotent. Its image
$z'_i$ in $R'$ therefore has positive valuation. The series
$\phi_i$ converge at $(z'_1,\ldots,z'_m)$, and the relations
$y_i=\phi_i$ give $\phi_i(z'_1,\ldots,z'_m)=b_i$, as required.
\end{proof}

For the bounded jumps $(a,M)$ treated later, we use $f_a^{(0)}$
when $M=pa$ and $f_a^{(1)}$ when $M>pa$, with the change of
variable above if $c\ne a$. These choices give $(p-1)a$ and
$(p-3)a$ new branch points, respectively, with nonzero finite
reduction in the coordinate $t$. In the second case the remaining
new branch points will be chosen separately.
Lemma~\ref{lem-finitecartier} gives the finiteness needed for the
later correction equations, and Lemma~\ref{lem-formal} allows
those equations to be solved after a finite extension of the
coefficient ring. 

\section{Ramification on a disk}\label{sec-analytic}

Let $k$ be an algebraically closed field of odd characteristic $p$.
Let $1\leq a<p$ be odd, and fix a lift $\chi_1$ of the degree-$p$
cover $w^p-w=x^{-a}$ together with its reflection. We seek a cyclic
cover of degree $p^2$ whose degree-$p$ quotient is $\chi_1$ and
whose reduction has the prescribed Witt coordinates $(x^{-a},f_2)$.
There are two parts to this problem. We must first choose the new
branch points so that the cover reduces to a smooth cover with the
required jumps $(a,M)$. We must then recover the coefficients of
$f_2$, since the jumps alone do not determine the cover.

In this section we give an equation in which the first quotient is
fixed and the new branch points can be varied. We explain how to test
a choice of branch points by computing a function
$\delta_\chi(r)$ on disks centered at the origin. Once the branch
points have been counted, the condition for smooth reduction with
jumps $(a,M)$ is simply that $\delta_\chi(r)=Mr$ for all sufficiently
small $r\geq0$. The last lemma gives the coefficient comparison that
will be needed to recover $f_2$.

Let $R$ be a complete discrete valuation ring of mixed characteristic
$(0,p)$ with residue field $k$, and let $L$ be its fraction field.
Write $v$ for the valuation normalized by $v(p)=1$. We allow finite
extensions of $L$, and continue to write $R$ for the valuation ring.
Let $T$ be the variable on the base of the lifted cover. Thus the
completed ring at the origin is $R[[T]]$, and reduction identifies
$R[[T]]/\mathfrak m_R R[[T]]$ with $k[[x]]$ by sending $T$ to $x$.
We choose $T$ so that the reflection sends $T$ to $-T$. After a
finite extension, $L$ contains a primitive $p^2$-th root of unity
$\zeta_{p^2}$. Set $\zeta_p=\zeta_{p^2}^p$ and $K=L(T)$.

We first recall how to write the cyclic covers of the $T$-line.
For $i=1,2$ and $A\in K^\times$, choose a root $Y$ of
$Y^{p^i}=A$ in a fixed algebraic closure of $K$. For
$\sigma\in\operatorname{Gal}(\overline K/K)$, define
$\kappa_i(A)(\sigma)\in\mathbb Z/p^i\mathbb Z$ by
\[
 \sigma(Y)=\zeta_{p^i}^{\,\kappa_i(A)(\sigma)}Y.
\]
The choice of root does not matter, since all the roots of unity
in this formula belong to $K$. This also shows that $\kappa_i(A)$
is a homomorphism. Its kernel fixes $K(Y)$, and the order of its
image is $[K(Y):K]$. We use additive notation for these
homomorphisms.

Multiplying roots gives
$\kappa_i(AB)=\kappa_i(A)+\kappa_i(B)$. Moreover,
$\kappa_i(A)=\kappa_i(B)$ if and only if $A/B$ is a $p^i$-th
power in $K$. Indeed, if the characters are equal, the quotient
of the chosen roots is fixed by every $\sigma$ and hence belongs
to $K$. The converse follows by multiplying a chosen root by
an element of $K$.

Every cyclic extension of degree $p^i$ occurs in this way. To see
this, let $\sigma$ generate its Galois group. As a $K$-linear map
on the extension field, $\sigma$ is diagonalizable, since
$X^{p^i}-1$ splits over $K$ with distinct roots. Some eigenvalue
has order $p^i$, since otherwise $\sigma^{p^{i-1}}$ would be the
identity. Raising a corresponding eigenvector to a suitable power
gives $Y$ with $\sigma(Y)=\zeta_{p^i}Y$. Then $Y^{p^i}\in K$,
and $Y$ has $p^i$ distinct conjugates, so it generates the
extension.

If $Y^{p^2}=A$, then $(Y^p)^p=A$, and an automorphism of value $j$
on $Y$ has value $j\bmod p$ on $Y^p$. Thus the character of the
degree-$p$ quotient is $\kappa_1(A)$. Identifying
$\mathbb Z/p\mathbb Z$ with $p\mathbb Z/p^2\mathbb Z$ by
$j\mapsto pj$, we write this as $p\kappa_2(A)=\kappa_1(A)$.

We shall also need to read the branch points directly from $A$.
Let $z$ be a geometric point of the $T$-line and set
$m=\operatorname{ord}_z(A)$, where a pole has negative order.
Then the ramification index at $z$ is
\[
 e_z=\frac{p^i}{\gcd(p^i,m)}.
\]
Indeed, in a local parameter $S$ at $z$ we may write $A=S^mu$,
where $u$ is a unit. Over $\overline L[[S]]$ the unit $u$ has a
$p^i$-th root by Hensel's lemma. It therefore suffices to consider
$Y^{p^i}=S^m$. If $d=\gcd(p^i,m)$, each resulting local field is
$\overline L((S^{1/(p^i/d)}))$, since $m/d$ and $p^i/d$ are
relatively prime. In particular, for a cover of degree $p^2$,
\[
 e_z=
 \begin{cases}
 p^2 & \text{if }p\nmid m,\\
 p   & \text{if }p\mid m\text{ and }p^2\nmid m,\\
 1   & \text{if }p^2\mid m.
 \end{cases}
\]
The middle case is the one that matters for our construction.
Such a point is branched in the degree-$p^2$ cover and unbranched
in its degree-$p$ quotient. We call it a new branch point.

We now choose an equation for the fixed quotient $\chi_1$.
Write $\chi_1=\kappa_1(A)$ and let $\tau$ act on $K$ by
$(\tau A)(T)=A(-T)$. Since the reflection acts on the cyclic
group by inversion, $\kappa_1(\tau A)=-\kappa_1(A)$. Consequently
$A(T)A(-T)$ is a $p$-th power. Taking its order at $z$ gives
\[
 \operatorname{ord}_{-z}(A)\equiv
 -\operatorname{ord}_z(A)\pmod p.
\]
At $0$ and $\infty$ this congruence says that twice the order is
zero modulo $p$. Since $p$ is odd, neither point is branched.
The branch points therefore occur in distinct pairs $\{z,-z\}$.
Choose one point from each pair and an integer $m_z$ congruent
to $\operatorname{ord}_z(A)$ modulo $p$. Set
\begin{equation}\label{eq-oldrep}
 A_0(T)=\prod_{\{z,-z\}}
        \left(\frac{T-z}{T+z}\right)^{m_z}.
\end{equation}
Then $A_0(-T)=A_0(T)^{-1}$. At every point, the order of $A/A_0$
is divisible by $p$. After splitting its zeros and poles, we may
factor its numerator and denominator into linear factors. This gives
$A/A_0=cB^p$ with $c\in L^\times$ and $B\in K^\times$. Adjoining
a $p$-th root of $c$ gives $\kappa_1(A_0)=\chi_1$.

We can now describe the freedom in the cover of degree $p^2$.
An equation $Y^{p^2}=C(T)$ has the prescribed quotient character
if and only if $\kappa_1(C)=\kappa_1(A_0)$, or equivalently
$C=A_0B^p$ for some $B\in K^\times$. If the reflection acts by
inversion on the full character, we may choose the equation in
the form
\begin{equation}\label{eq-kummer}
 Y^{p^2}=A_0(T)G(T)^p
 \quad\text{with}\quad G(-T)=G(T)^{-1}.
\end{equation}
Here is why the last equality can be imposed. Since $C\,\tau C$
is a $p^2$-th power and $A_0\,\tau A_0=1$, the function
$(B\,\tau B)^p$ is a $p^2$-th power. Thus $B\,\tau B$ is a
$p$-th power, since every $p$-th root of unity in $K$ is itself
a $p$-th power in $K$. Modulo $p$-th powers we therefore have
$\tau B\equiv B^{-1}$. Set
\[
 G=\left(\frac{B}{\tau B}\right)^{(p+1)/2}.
\]
Then $\tau G=G^{-1}$, and $G\equiv B^{p+1}\equiv B$ modulo
$p$-th powers. Hence $A_0G^p$ and $A_0B^p$ define the same
character of order $p^2$.

Conversely, the tower of function fields defined by
\eqref{eq-kummer} is
\[\begin{tikzcd}[column sep=large]
 K \arrow[r,hook,"p"] &
 K(W) \arrow[r,hook,"p"] &
 K(Y),
\end{tikzcd}\]
where
 \[W=\frac{Y^p}{G(T)},\qquad W^p=A_0(T).
\]
In particular the quotient is independent of $G$. The reflection
sends $(T,Y)$ to $(-T,Y^{-1})$ and induces
$(T,W)\mapsto(-T,W^{-1})$. Two choices $G$ and $G'$ define the
same character precisely when $G'/G$ is a $p$-th power. The
formula classifies covers with the fixed quotient character over
$L$; the condition on their reduction is still to be imposed.

For example, take $b\ne0$ away from the zeros and poles of $A_0$,
and set $G(T)=(T-b)/(T+b)$. The orders of $A_0G^p$ at $b$ and
$-b$ are $p$ and $-p$, respectively. Thus we have added two branch
points of index $p$ without changing the quotient. Varying $b$
moves this pair. Which choices can give the required reduction is
the question addressed by the depth function.

For rational $r\geq0$, let
\[
 D[r]=\{z\in\overline L\mid v(z)\geq r\}
 \quad\text{and}\quad
 D(r)=\{z\in\overline L\mid v(z)>r\},
\]
where $v(0)=+\infty$. With $|z|=p^{-v(z)}$, these are the closed
and open disks of radius $p^{-r}$ centered at zero. Increasing $r$
makes the disk smaller. In particular, $D(0)$ consists of the
points reducing to $x=0$.

Set $\lambda_p=\zeta_p-1$ and $\alpha=p/(p-1)$. Choose compatible
roots of $\lambda_p$ in $\overline L$. This gives elements
$\varpi_s$ for $s\in\mathbb Q$ such that $v(\varpi_s)=s$,
$\varpi_{s+t}=\varpi_s\varpi_t$, and
$\varpi_{1/(p-1)}=\lambda_p$. Only finitely many of them are needed
at any one time. In particular $\varpi_\alpha=\lambda_p^p$.

To study $D[r]$, set $T_r=T/\varpi_r$ and let $t$ be the reduction
of $T_r$. A point $z\in D[r]$ then gives the point
$t=\overline{z/\varpi_r}$ on the affine $t$-line. This point is
zero if $v(z)>r$ and nonzero if $v(z)=r$. More generally, two
points $z,z'\in D[r]$ have the same reduction in $t$ if and only
if $v(z-z')>r$. In particular, two points of the same valuation
need not have the same reduction.

For $P=\sum_j a_jT_r^j$, set $v_r(P)=\min_j v(a_j)$.
This rule is multiplicative. Indeed, after scaling two polynomials
so that the minimum is zero for each, their reductions are nonzero
polynomials in $k[t]$, and their product is nonzero. We may therefore
set $v_r(P/Q)=v_r(P)-v_r(Q)$. This is the Gauss valuation on $K$.
Its residue field is $k(t)$. For example, if $v_r(P/Q)=0$, scale
$P$ and $Q$ by the same constant so that both have valuation zero.
The reduction of $P/Q$ is then $\bar P(t)/\bar Q(t)$.
Substitution of $T=\varpi_rT_r$ gives
\begin{equation}\label{eq-gauss}
 v_r\left(\sum_{j\geq0}c_jT^{-j}\right)
 =\inf_{j\geq0}\{v(c_j)-jr\}
\end{equation}
whenever the series converges for $v_r$. For a rational function
with a pole on $v(T)=r$, we use a quotient of polynomials directly.

We next describe the two quantities attached to a character $\chi$
at $v_r$. The first is a number $\delta_\chi(r)$, its depth. The
second, when the depth is positive, is a rational differential
$\omega_\chi(r)$ in the variable $t$. We use the normalization of
\cite[Section~6.3]{Obus}. At a fixed rational $r$ we may extend
the constants so that the extension of value groups has index one.
Depth zero then means that the character is unramified for $v_r$
and comes by reduction from a character over $k(t)$.

For precision, the positive depths used below can be defined as
follows. Complete $K$ at $v_r$. In the case where the residue
extension of the cyclic cover is purely inseparable of the full
degree, let $E$ be the resulting extension of complete fields,
$\mathcal O_E$ its valuation ring and $J$ its Galois group.
Keep the normalization $v(p)=1$. For $s\geq0$ set
\[
 \begin{split} &J_s=\{\sigma\in J\mid
          v(\sigma(u)-u)\geq s\ \text{for every }u\in\mathcal O_E\},\\
 & h(s)=\int_0^s |J_b|\,db .
 \end{split}
\]
The depth is the largest $h(s)$ for which $\chi$ is nonzero on
$J_s$. This is the upper ramification break with the valuation
normalization used here \cite[Section~3.1 and Definition~3.3]{Wewers}.
The differential is normalized using the same elements $\varpi_s$,
as in \cite[Definition~3.9]{Wewers}. For the equations in this paper
it is computed by the following formulas.

Let $F\in K^\times$ have $v_r(F)=0$, and suppose that its reduction
$\bar F$ is not a $p$-th power in $k(t)$. If $0<b<\alpha$, then
\begin{equation}\label{eq-kummerdepth}
 \delta_{\kappa_1(1+\varpi_bF)}(r)=\alpha-b
 \quad\text{and}\quad
 \omega_{\kappa_1(1+\varpi_bF)}(r)=d\bar F .
\end{equation}
Here $d\bar F=\bar F'(t)\,dt$. We explain the depth calculation.
Set $\theta=\varpi_{b/p}$ and substitute $Y=1+\theta V$ in
$Y^p=1+\varpi_bF$. Dividing by $\theta^p$ gives an integral
equation for $V$ whose reduction is $\bar V^p=\bar F$. Indeed,
for $1\leq j<p$ the coefficient of $V^j$ has valuation
$1-(p-j)b/p>0$. Since $\bar F$ is not a $p$-th power,
$1,\bar V,\ldots,\bar V^{p-1}$ are linearly independent over
$k(t)$. It follows that $1,V,\ldots,V^{p-1}$ are an integral
basis. To see this, scale a linear combination so that the least
coefficient valuation is zero. Its reduction is nonzero by linear
independence, so its valuation is zero as well.

For each nonidentity automorphism $\sigma$ of this degree-$p$
extension, $v(\sigma(Y)-Y)=1/(p-1)$, since $Y$ is a unit and
$\sigma(Y)/Y$ is a primitive $p$-th root of unity. Hence
\[
 v(\sigma(V)-V)=\frac{1}{p-1}-\frac{b}{p}.
\]
Every difference $\sigma(u)-u$ with $u\in\mathcal O_E$ is divisible
by $\sigma(V)-V$, by the integral basis just obtained. Thus $J_s$
has order $p$ up to the displayed value of $s$ and is trivial
thereafter. Multiplying that value by $p$ gives $\alpha-b$.
The differential formula in \eqref{eq-kummerdepth} is
\cite[Proposition~5.3(ii)]{Wewers}. For a unit equation $Y^p=F$
with the same hypothesis on $g=\bar F$, the corresponding formulas
are $\delta=\alpha$ and $\omega=dg/g$
\cite[Proposition~5.3(i)]{Wewers}.

At $b=\alpha$ a different calculation applies. Substituting
$Y=1+\lambda_p y$ and dividing by $\lambda_p^p$ gives the reduction
\[
 y^p-y=\bar F.
\]
Indeed, $p/\lambda_p^{p-1}$ reduces to $-1$, by
$(1+\lambda_p)^p=1$, and the coefficients of $y^j$ for
$1<j<p$ have positive valuation. The derivative of $y^p-y-\bar F$
is $-1$, so this equation is separable and the depth is zero. It
may even split if $\bar F=h^p-h$ for some $h\in k(t)$. This
calculation concerns the valuation $v_r$ only. Smooth reduction of
the whole cover requires the additional test given below.

The formulas are particularly useful when $G$ is changed. Suppose
that $\chi=\psi+\psi'$ and fix $r$. If the two depths are unequal,
$\chi$ has the larger depth and the differential of that summand.
If both depths are equal to $\delta>0$, add their differentials.
A nonzero sum is the differential of $\chi$, whose depth is then
$\delta$. If the sum is zero, the depth of $\chi$ is strictly less
than $\delta$. If both depths are zero, the reductions add as
characters over $k(t)$. These are the addition rules of
\cite[Proposition~6.20]{Obus}.

As $r$ varies, $\delta_\chi(r)$ is continuous and piecewise linear,
and changes of slope occur only at rational values of $r$. At a
rational $r>0$ of positive depth, its slopes are
\begin{equation}\label{eq-slopes}
 \delta_\chi'(r^-)=\operatorname{ord}_\infty\omega_\chi(r)+1
 \quad\text{and}\quad
 \delta_\chi'(r^+)=-\operatorname{ord}_0\omega_\chi(r)-1
\end{equation}
by \cite[Proposition~6.17]{Obus}. The order of a differential at a
point is the order of its coefficient when written in a local
parameter at that point. For example, a term $ct^{-j-1}dt$ with
$c\ne0$ has order $-j-1$ at zero and order $j-1$ at infinity.
The latter follows by setting $u=1/t$, which gives
$ct^{-j-1}dt=-cu^{j-1}du$.

There is also a bound at each finite point $b$ of the residue line.
Let $n_b$ be the number of branch points $z\in D[r]$ satisfying
$\overline{z/\varpi_r}=b$. Then
$\operatorname{ord}_b\omega_\chi(r)\geq-n_b$
by \cite[Proposition~6.19]{Obus}. In particular, a residue class
containing no branch point gives no pole. If every branch point
lies in $D(r)$, the only possible finite pole is therefore at zero.
In our covers all branch points are finite and nonzero. To check
regularity at infinity, choose $d>r$ larger than all their valuations
and set $V=\varpi_d/T$. At radius $d-r$ the new residue variable
is $1/t$, since $V/\varpi_{d-r}=\varpi_r/T$. Its zero class would
contain a branch point $z$ only if $v(z)<r$. There is no such point.

It follows that if all branch points lie in $D(r)$ and
$\omega_\chi(r)$ is exact, then it has a primitive in
$t^{-1}k[t^{-1}]$. Indeed, regularity away from zero and at infinity
gives a finite expansion $\omega_\chi(r)=\sum_{j>0}c_jt^{-j-1}dt$.
The derivative of a Laurent series in characteristic $p$ has zero
coefficient at $t^{-j-1}dt$ whenever $p\mid j$. Thus $c_j=0$
for these $j$, and a primitive is
$-\sum_{p\nmid j}(c_j/j)t^{-j}$.

We can now state the test at the original radius $r=0$. Suppose
that the degree-$p$ quotient has good reduction with jump $a$,
that every geometric branch point of the cyclic cover lies in
$D(0)$, and that there are $a+1$ points of index $p^2$ and $M-a$
points of index $p$. Good reduction means that the cover extends
with smooth special fibre and that its reduction retains the full
degree. Under the stated branch counts,
\cite[Proposition~6.18]{Obus} gives
\[
 \text{good reduction with jumps }(a,M)
 \quad\Longleftrightarrow\quad
 \delta_\chi(0)=0\ \text{ and }\ \delta_\chi'(0^+)=M.
\]
Since the depth is piecewise linear, the two conditions on the
right are equivalent to $\delta_\chi(r)=Mr$ for all sufficiently
small $r\geq0$. Notice that the hypothesis on a branch point is
$v(z)>0$. This ensures that it reduces to $x=0$.

We return to the fixed quotient $\chi_1$. Set $r_c=1/((p-1)a)$.
The lift in \cite[Lemma~4.2]{Obus} has all its $a+1$ branch points
in $D(r_c)$. Since it has good reduction,
\cite[Lemma~6.1]{OW} gives
\begin{equation}\label{eq-firstdepth}
 \delta_{\chi_1}(r)=ar \quad\text{for }0\leq r\leq r_c .
\end{equation}
For rational $0<r\leq r_c$, the same lemma gives
$\omega_{\chi_1}(r)=c_r\,dt/t^{a+1}$ for some $c_r\in k^\times$.
The inverse image of $D[r]$ in the degree-$p$ cover is a closed
disk. Indeed, $D(r)$ contains a branch point of index $p$, so its
inverse image is connected by \cite[Lemma~2.16]{Kinks}. The left
slope is $a$, one less than the number of branch points in $D[r]$.
The assertion now follows from \cite[Corollary~2.18]{Kinks}.

Let $\chi_2$ be a character of order $p^2$ with $p\chi_2=\chi_1$.
Write $\delta_i=\delta_{\chi_i}(r)$ and, at positive depth,
$\omega_i=\omega_{\chi_i}(r)$. The reason for singling out $r_c$
is that $\delta_1$ reaches $1/(p-1)$ there. Before it reaches
this value, if $0<\delta_1<1/(p-1)$ and $\delta_2<\alpha$,
\cite[Theorem~4.3(iii)]{Wewers} gives $\delta_2\geq p\delta_1$
and
\begin{equation}\label{eq-cartierdepth}
 \Cart(\omega_2)=
 \begin{cases}
 \omega_1 & \text{if }\delta_2=p\delta_1,\\
 0 & \text{if }\delta_2>p\delta_1.
 \end{cases}
\end{equation}
Here $\Cart$ is the Cartier operator. Explicitly, writing
$H(t)=\sum_{j=0}^{p-1}h_j(t)^pt^j$ gives
$\Cart(H(t)\,dt)=h_{p-1}(t)\,dt$. Thus in the second case
$\omega_2$ is exact, since the terms with $j<p-1$ can each be
integrated. At $r=r_c$ the same theorem forces $\delta_2=\alpha$
and gives
\begin{equation}\label{eq-criticalcartier}
 \Cart(\omega_2)-\omega_2
 =\omega_1=c_{r_c}\,\frac{dt}{t^{a+1}}.
\end{equation}
This explains where the equation of the preceding section enters.
For a solution of the form $dt/(t^{a+1}f(t))$, the nonzero roots
of $f$ specify the reductions $\overline{z/\varpi_{r_c}}$ of the
new branch points with $v(z)=r_c$. They specify these residues,
whereas the points $z$ themselves must still be chosen. To pass
from the coefficient $a$ in Lemma~\ref{lem-data} to $c_{r_c}$,
choose $\lambda\in k^\times$ with $\lambda^a=c_{r_c}/a$ and replace
$f(t)$ by $\lambda^{-a}f(t/\lambda)$. This multiplies the roots by
$\lambda$ while leaving $T$ and its reduction $x$ fixed.

The remaining choice is made by varying the points in a family.
The following theorem explains why one can choose a member for
which the desired slope persists as far as possible. Its conclusion
is that a minimum is attained. Showing that a positive minimum
can be improved is a separate step in the lifting proof.

\begin{theorem}[{\cite[Corollaries~2.20(ii)--(iv) and~5.3(ii)]{Kinks}}]
\label{thm-minimum}
Let $\mathcal A$ be a nonempty rigid analytic parameter space with
a finite covering by affinoid subdomains whose pairwise intersections
also admit finite affinoid coverings. Let
$\mathcal Y\to\mathcal Y_1\to\mathbb P^1_{\mathcal A}$ be a family
of cyclic covers of degree $p^2$ with their degree-$p$ quotients.
Assume that both maps are finite flat Galois covers of smooth proper
curves over $\mathcal A$, and that their branch points form finite
\'etale families of constant degree. Use the same coordinate $T$ in
each fibre.

Suppose there are rational numbers $s_1,s_2>0$ for which the
following conditions hold in every fibre.
\begin{enumerate}
\item The inverse image of $D(r)$ in the full cover is connected for
$0<r<s_1$.
\item There are exactly $d$ branch points of the full cover in
$D(0)$, all of which belong to $D[s_2]$.
\item Set $r_0=\min(s_1,s_2)$. For rational $0<r\leq r_0$, the
inverse image of $D[r]$ in the degree-$p$ quotient is a closed disk.
At each such $r$, after extending constants, $v_r$ has a unique
extension to the full cover and its residue field is purely
inseparable of degree $p^2$ over $k(t)$.
\end{enumerate}
Let $\chi$ be a faithful character of a fibre. For rational
$0<r\leq r_0$ one has $\delta_\chi'(r^-)\leq d-1$, with equality
if and only if the inverse image of $D[r]$ in the full cover is
a closed disk. Let $\lambda(\chi)$ be the largest such $r$ for
which $\delta_\chi'(r^-)<d-1$, taking $\lambda(\chi)=0$ if
there is none. Then $\lambda(\chi)$ attains its minimum on a
nonempty subdomain of $\mathcal A$.
\end{theorem}

In this statement an affinoid subdomain is an analytic piece
described by convergent equations on a closed polydisk. The
finiteness hypothesis on their covering is often abbreviated
\emph{qcqs}. The condition on the branch points means that locally,
after a finite \'etale base change, they can be listed as disjoint
sections. Thus the count $d$ counts distinct points in every fibre.
A residue field extension is purely inseparable if every element
in it has a $p$-power in the base field.

We finish with the coefficient comparison. In the next lemma $F_2$
is the function obtained directly from the branch points.
Dividing it by a $p$-th power does not change its character and
allows us to remove certain terms whose degrees are divisible
by $p$. The point of the lemma is that this operation preserves
the coefficient in degree $m$ up to terms of valuation greater
than $\beta$, and has the same property in every higher degree.

\begin{lemma}\label{an-cutoff}
Fix a rational $r\geq0$, set $U=T_r$, and let
$0<\beta\leq\alpha$ be rational. Let $m>0$ be prime to $p$.
Suppose that $F_2=1+\sum_{i>0}b_iU^{-i}$ is a rational function
whose expansion converges on $|U|\geq1$. Assume that $v(b_i)>0$
for all $i>0$, that $v(b_m)\geq\beta$, and that $v(b_i)>\beta$
for $i>m$.

If $\delta_{\kappa_1(F_2)}(r)\leq\alpha-\beta$, then, after a
finite extension of $L$, there exists a rational function $H$ such
that $F_1=F_2/H^p=1+\sum_{i>0}a_iU^{-i}$ converges on the same
domain and satisfies
\[
 v(a_i)\geq\beta\text{ for every }i>0
 \quad\text{and}\quad
 v(a_i)>\beta\text{ if }p\mid i.
\]
For every such $F_1$, one has $v(a_i-b_i)>\beta$ for all $i\geq m$.
\end{lemma}

\begin{proof}
Set $X=U^{-1}$ and write $B(X)=1+\sum_{i>0}b_iX^i$. Choose $n$
with $pn\geq m$. We first construct a polynomial
$H=1+\sum_{j=1}^n h_jX^j$ such that $v(h_j)>0$ and
\[
 [X^{pj}]H^p=b_{pj}\quad\text{for }1\leq j\leq n.
\]
Here $[X^l]$ denotes the coefficient of $X^l$. To see existence,
homogenize these $n$ equations by introducing $h_0$. They become
\[
 [X^{pj}](h_0+h_1X+\cdots+h_nX^n)^p=b_{pj}h_0^p
 \quad\text{for }1\leq j\leq n.
\]
Over $\overline L$, $n$ homogeneous equations of positive degree
in $n+1$ variables have a common nonzero solution. This follows
from the height theorem, since the ideal of the origin has height
$n+1$. The coordinates of a solution lie in a finite extension of
$L$. Scale them so that all are integral and at least one is a
unit. On reduction the equations give $\bar h_j^p=0$ for $j>0$,
because $v(b_{pj})>0$. Hence $h_0$ is a unit. Dividing the
coordinates by $h_0$ gives the required $H$. In particular
$\bar H=1$, so $H$ has no zero on $|X|\leq1$.

Set $D=B-H^p$. Every coefficient of $D$ in positive degree
divisible by $p$ has valuation greater than $\beta$. In degrees
$p,2p,\ldots,pn$ these coefficients vanish by construction. Beyond
degree $pn$ they are the coefficients of $B$, which have valuation
greater than $\beta$ since $pn\geq m$.

We claim that $v_r(D)\geq\beta$. Otherwise set $s=v_r(D)<\beta$.
The constant term of $D$ is zero and all other coefficients have
positive valuation, so $s>0$. The reduction of $D/\varpi_s$ is
a nonzero polynomial supported in positive degrees not divisible
by $p$. It cannot be a $p$-th power in $k(t)$. Since $\bar H=1$,
the reduction of $(B/H^p-1)/\varpi_s$ is the same polynomial.
Equation~\eqref{eq-kummerdepth} now gives
\[
 \delta_{\kappa_1(F_2)}(r)
 =\delta_{\kappa_1(B/H^p)}(r)
 =\alpha-s>\alpha-\beta,
\]
a contradiction. This proves the claim. Since
$F_1-1=D/H^p$ and $H$ reduces to $1$, the coefficients of $F_1-1$
have valuation at least $\beta$, with strict inequality in degrees
divisible by $p$.

It remains to compare an arbitrary such $F_1$ with $F_2$. This is
the comparison of \cite[Lemma~7.63]{Obus}; we include the argument
to account for the endpoint $\beta=\alpha$. Write $F_2=F_1J^p$.
Both $F_1$ and $F_2$ are units on $|X|\leq1$ reducing to $1$.
Hence $J$ has no zero or pole there, and its reduction satisfies
$\bar J^p=1$, so $\bar J=1$. Multiplying $J$ by a $p$-th root of
unity, we may arrange that $J(0)=1$. Write
$J=1+\sum_{i>0}u_iX^i$. Then $v(u_i)>0$ and $v(u_i)\to+\infty$.

Let $j$ be the largest index for which $v(u_j)\leq\beta/p$, and
set $j=0$ if there is none. Write $J=P+Q$, where $P$ consists of
the terms through degree $j$. Every coefficient of $Q$ has
valuation greater than $\beta/p$. In
\[
 J^p-P^p=Q^p+\sum_{\ell=1}^{p-1}\binom p\ell P^{p-\ell}Q^\ell
\]
all coefficients have valuation greater than $\beta$. For $Q^p$
this follows by multiplying $p$ coefficients of $Q$. For a mixed
term it follows from $p\mid\binom p\ell$ and
$1+\beta/p\geq\beta$, the latter inequality being equivalent to
$\beta\leq p/(p-1)$.

Consequently, if $j>0$, the coefficient of $X^{pj}$ in $J^p-1$
has valuation $p\,v(u_j)\leq\beta$, and all coefficients in higher
degrees have valuation greater than $\beta$. Also every coefficient
of $(F_1-1)(J^p-1)$ has valuation greater than $\beta$. If $pj>m$,
then in the identity
\[
 F_2-1=(F_1-1)+(J^p-1)+(F_1-1)(J^p-1)
\]
the first and third terms have coefficient of valuation greater
than $\beta$ in degree $pj$, whereas the second does not. This
contradicts $v(b_{pj})>\beta$. Equality $pj=m$ is impossible
since $p\nmid m$. Thus $pj<m$, including when $j=0$. It follows
that every coefficient of $J^p-1$ in degree at least $m$ has
valuation greater than $\beta$. The identity above now gives
$v(b_i-a_i)>\beta$ for every $i\geq m$, as claimed.
\end{proof}

\section{The new branch points}\label{sec-normalized}

Let $k$ be an algebraically closed field of odd characteristic $p$.
We consider a local $D_{p^2}$-action whose cyclic upper jumps are
$(a,M)$, where $1\leq a<p$ and $M=pa+2e$ with $0\leq e<p$.
Recall that $a$ is odd. We have chosen the first Witt coordinate to
be $x^{-a}$, so the quotient of degree $p$ is given by
$w^p-w=x^{-a}$. We fix a lift of this quotient. The question is how
to extend it to degree $p^2$ so that its reduction has the given
second Witt coordinate $f_2$.

Write the fixed lift as $W^p=A_0(T)$, where $T$ reduces to $x$ and
$A_0(-T)=A_0(T)^{-1}$. We seek a rational function $G(T)$ such that
\[
 Y^{p^2}=A_0(T)G(T)^p
 \quad\text{and}\quad
 G(-T)=G(T)^{-1}.
\]
Indeed, on setting $W=Y^p/G(T)$ we recover $W^p=A_0(T)$, and the
reflection acts by $(T,Y)\mapsto(-T,Y^{-1})$. Thus we may vary $G$
without changing either the quotient or its reflection.

In this section we construct $G$ with the required number of branch
points and determine which of its coefficients can be varied.
Section~\ref{sec-corr-formal} gives the equations for these
variations. In Section~\ref{sec-boundedproof} we use them to obtain
good reduction and to prescribe $f_2$. This proves the following
proposition.

\begin{proposition}\label{prop-normalized}
Let $p$ be an odd prime, let $1\leq a<p$ be odd, and let
$0\leq e<p$. Every local $D_{p^2}$-action with cyclic upper jumps
$(a,pa+2e)$ lifts to characteristic zero.
\end{proposition}

Let $R$ be the valuation ring of the fixed lift and normalize its
valuation $v$ by $v(p)=1$. We allow finite extensions of its fraction
field. All branch points counted below are points on the $T$-line.
For good reduction with jumps $(a,M)$, the required counts are
\[
 \begin{array}{c|c|c}
 \text{ramification index} & \text{number of points}
                           & \text{source of the points}\\ \hline
 p^2 & a+1 & \text{branch points of } W^p=A_0(T)\\
 p   & M-a & \text{new branch points}
 \end{array}
\]
and all these points must reduce to $x=0$. Once the counts hold,
the criterion of \cite[Proposition~6.18]{Obus} takes the form
\[
 \text{good reduction with jumps }(a,M)
 \ \Longleftrightarrow\
 \delta(0)=0\ \text{ and }\ \delta'(0^+)=M,
\]
where $\delta(r)$ is the depth of the degree-$p^2$ cover on the disk
$D[r]=\{v(T)\geq r\}$. Thus the count tells us how many points to
choose. It does not by itself give good reduction.

Set $r_c=1/((p-1)a)$. The $a+1$ branch points of the fixed quotient
satisfy $v(T)>r_c$. If $e=0$, then $M-a=(p-1)a$, and the first
polynomial of Lemma~\ref{lem-data} has precisely this degree. Its
roots supply all the new points at $v(T)=r_c$. No second radius is
needed in this case.

Suppose for now that $e>0$. Set
\begin{equation}\label{eq-parameters}
 \begin{aligned}
 N&=(p-3)a, & n&=N/2, & K&=a+N=(p-2)a,\\
 B&=a+e,    & E&=B-1, & M&=K+2B.
 \end{aligned}
\end{equation}
The second polynomial of Lemma~\ref{lem-data} has degree $N$. Its
roots give $N$ new points at $v(T)=r_c$. Thus $K+1=(a+1)+N$ is
the number of points accounted for so far. The number still to be
chosen is
\[
 (M-a)-N=M-K=2B.
\]
If $p=3$, then $N=0$ and all the new points belong to this last
group.

We choose the $2B$ remaining points to have a common valuation
$r_h<r_c$. To determine $r_h$, set $\alpha=p/(p-1)$. At $r_c$ the
fixed quotient has depth $ar_c=1/(p-1)$, so the cover of degree
$p^2$ has depth $\alpha$ by \cite[Theorem~4.3]{Wewers}. If the
cover has good reduction, the slope on an interval containing no
branch valuation is the number of branch points in the disk minus
one \cite[Propositions~6.17 and~6.19]{Obus}. For $0<r<r_h$ the disk
contains all $M+1$ points. For $r_h<r<r_c$ it contains $K+1$.
The required depth is therefore
\[
 \delta(r)=
 \begin{cases}
 Mr & 0\leq r\leq r_h,\\
 \alpha-K(r_c-r) & r_h\leq r\leq r_c.
 \end{cases}
\]
The two expressions must agree at $r_h$. Using $M-K=2B$ and
$\alpha-Kr_c=2/(p-1)$, we obtain
\begin{equation}\label{eq-parameters2}
 r_h=\frac{\alpha-Kr_c}{M-K}=\frac{1}{(p-1)B},
 \quad\text{and}\quad
 r_c=\frac{1}{(p-1)a}.
\end{equation}
Since $B=a+e>a$, we have $0<r_h<r_c$. Set $\Delta=r_c-r_h$ and
$s_h=K\Delta$. We shall use the identities
\begin{equation}\label{eq-corr-budget}
 K+2B=M
 \quad\text{and}\quad
 s_h+Mr_h=\alpha.
\end{equation}
For example, take $p=5$ and $a=e=1$. The quotient supplies two
points of index $25$. We need six more points, all of index $5$.
Two have valuation $1/4$ and four have valuation $1/8$.

We now construct the $2B$ points of valuation $r_h$. The following
observation allows us to produce several points at once. A pole of
order $J$ gives one branch point at the pole itself and $J$ more at
nearby zeros.

\begin{lemma}\label{lem-auxiliary-branches}
Let $1\leq J<p$, let $\xi,a_1,\ldots,a_J\in R$ with
$a_J\in R^\times$, and let $v(\varepsilon)>0$ with
$\varepsilon\ne0$. Set
\[
 A(t)=\sum_{j=1}^J\frac{a_j}{(t-\xi)^j}.
\]
The cover $Y^p=1+\varepsilon A(t)$ has exactly $J+1$ geometric
branch points. One is $\xi$. The others are simple zeros
$t_1,\ldots,t_J$ satisfying
\[
 v(t_i-\xi)=\frac{v(\varepsilon)}{J}.
\]
For fixed $\varepsilon$, the zeros form a finite \'etale family
over the characteristic-zero parameter space as $\xi$ and the
$a_j$ vary subject to the conditions above.
\end{lemma}

\begin{proof}
Choose $\rho$ with $\rho^J=\varepsilon$ and write $t=\xi+\rho u$.
After multiplication by $u^J$, the equation
$1+\varepsilon A(t)=0$ becomes
\[
 u^J+\sum_{j=1}^J a_j\rho^{J-j}u^{J-j}=0.
\]
Its reduction is $u^J+\bar a_J$. Since $p\nmid J$ and
$\bar a_J\ne0$, this polynomial has $J$ distinct nonzero roots.
Each lifts by Hensel's lemma. Thus the zeros are simple and have
the form $\xi+\rho u_i$ with $v(u_i)=0$. The same calculation shows
that the discriminant of the polynomial in $u$ is a unit throughout
the parameter space, which gives the assertion about the family.

The function $1+\varepsilon A(t)$ has a pole of order $J$ at $\xi$
and takes the value $1$ at infinity. A cover $Y^p=H(t)$ is branched
at a zero or pole of $H$ precisely when its order is not divisible
by $p$. This gives the stated $J+1$ points.
\end{proof}

The same calculation applies to a sum of such functions $A$ provided
their centres $\xi$ have distinct reductions. Near each centre the
other summands are regular and integral, so they do not change the
reduced equation $u^J+\bar a_J=0$. The zeros found in this way
account for all the zeros, since their number is the sum of the pole
orders.

For each valuation $s$ used below, choose an element $\varpi_s$ with
$v(\varpi_s)=s$, compatibly as in the preceding section. At $r_h$
we use the variable $t=\overline{T/\varpi_{r_h}}$. By
\eqref{eq-slopes}, the required slopes $K$ on the right and $M$ on
the left amount to the following two conditions on the differential
$\eta$ at $r_h$,
\[
 \operatorname{ord}_0(\eta)=-K-1
 \quad\text{and}\quad
 \operatorname{ord}_\infty(\eta)=M-1.
\]
We construct $\eta$ as $dP=P'(t)\,dt$. The terms of $P$ with poles
away from zero will give the new branch points. We also need to know
which coefficients of these terms can be varied independently.

\begin{lemma}\label{lem-auxiliary}
For each $C\in k^\times$, there is an odd rational function $P(t)$
such that $dP$ has leading term $Ct^{-K-1}dt$ at zero and order
$M-1$ at infinity. The nonzero poles of $P$ form one or two pairs
$\{\xi,-\xi\}$. Each has order less than $p$, and together they give
$2B$ branch points by Lemma~\ref{lem-auxiliary-branches}.

Let $P_{\mathrm{aux}}$ be the sum of the principal parts of $P$ at
its nonzero poles. Here the principal part at $\xi$ is the sum of the
terms $(t-\xi)^{-j}$, with their coefficients, for $j>0$. Write
$[t^{-q}]P_{\mathrm{aux}}$ for the coefficient of $t^{-q}$ in its
expansion at infinity. Among the pole positions and the coefficients
of their principal parts, there are $E$ parameters for which the
Jacobian of the following $E$ coefficients is invertible,
\begin{equation}\label{eq-corr-aux-rows}
 [t^{-q}]P_{\mathrm{aux}}
 \quad\text{for}\quad
 q\in\{K,K+2,\ldots,M-2\}\setminus\{pa\}.
\end{equation}
\end{lemma}

\begin{proof}
We first construct $\eta$. If $B\leq p$, choose $b\in k^\times$
and set
\begin{equation}\label{eq-oneaux}
 \eta=\frac{C(-b)^B\,dt}{t^{K+1}(t^2-b)^B}.
\end{equation}
If $B>p$, set $b_2=B-p+1$. Then $2\leq b_2\leq a<p$. Choose
$c_2\in k^\times$ and set
\[
 b=\frac{ac_2}{a-b_2+1}.
\]
In particular, $b\ne0$ and $b\ne c_2$. Set
\begin{equation}\label{eq-twoaux}
 \eta=
 \frac{C(-b)^{p-1}(-c_2)^{b_2}\,dt}
 {t^{K+1}(t^2-b)^{p-1}(t^2-c_2)^{b_2}}.
\end{equation}
In both cases the numerator gives leading term $Ct^{-K-1}dt$
at zero. The denominator has degree $K+1+2B=M+1$. Substituting
$t=1/u$ and $dt=-u^{-2}du$ gives order $M-1$ at infinity.

We show that $\eta=dP$ for a rational function $P$. Write
$\eta=c\,dt/(t^{K+1}D(t^2))$, where $c=C D(0)$ and
\[
 D(z)=
 \begin{cases}
 (z-b)^B & B\leq p,\\
 (z-b)^{p-1}(z-c_2)^{b_2} & B>p.
 \end{cases}
\]
Since $K=pa-2a$, the substitution $z=t^2$ gives
\[
 \eta=\frac c2\,t^{-pa}
       \left.\frac{z^{a-1}}{D(z)}\,dz\right|_{z=t^2}.
\]
In characteristic $p$ we have $d(t^{-pa})=0$. It therefore suffices
to find $V(z)$ such that $dV=z^{a-1}dz/D(z)$.

If $B\leq p$, expand $z^{a-1}=(b+(z-b))^{a-1}$. After division
by $(z-b)^B$, the only terms are multiples of $(z-b)^{-j}$ with
$e+1\leq j\leq B$. Since $2\leq j\leq p$, each has the primitive
$(z-b)^{1-j}/(1-j)$.

Suppose $B>p$. The reason for the choice of $b$ is that the
coefficient of $(z-c_2)^{-1}$ must vanish. Indeed, write
$(z-b)^{1-p}=(z-b)/(z-b)^p$. In the expansion at $c_2$, the
factor $(z-b)^{-p}$ has no terms of degrees $1,\ldots,p-1$.
The required coefficient is therefore
\begin{equation}\label{eq-residue}
 \begin{aligned}
 \operatorname{Res}_{z=c_2}\frac{z^{a-1}\,dz}{D(z)}
 &=(c_2-b)^{-p}
   [u^{b_2-1}](c_2+u)^{a-1}(c_2+u-b)\\
 &=(c_2-b)^{-p}\binom{a-1}{b_2-1}c_2^{a-b_2}
   \left(\frac{ac_2}{a-b_2+1}-b\right)=0 .
 \end{aligned}
\end{equation}
The factor outside the parentheses is nonzero, so this equality
also determines $b$ from $c_2$. The differential is regular at
infinity. Its residue at $b$ is therefore zero as well, since the
sum of the residues is zero. All remaining partial fractions have
the form $h(z-\xi)^{-j}$ with $2\leq j<p$, and can be integrated
term by term.

In either case take $V$ to be the sum of the primitives just
constructed. It is regular at zero. Set
\[
 P(t)=\frac c2\,t^{-pa}\bigl(V(t^2)-V(0)\bigr).
\]
Then $P(-t)=-P(t)$ and $dP=\eta$. At each nonzero pole, integration
lowers the pole order by one. Lemma~\ref{lem-auxiliary-branches}
now gives the following counts,
\[
 \begin{array}{c|c|c|c}
 & \text{centres} & \text{pole order of }P
                  & \text{branch points per centre}\\ \hline
 B\leq p & \pm\sqrt b   & B-1   & B\\ \hline
 B>p     & \pm\sqrt b   & p-2   & p-1\\
         & \pm\sqrt{c_2}& b_2-1 & b_2
 \end{array}
\]
Thus there are $2B$ points in the first case and
$2(p-1)+2b_2=2B$ in the second.

It remains to prove the assertion about the coefficients in
\eqref{eq-corr-aux-rows}. At one centre in each pair, vary the
coefficients of the principal part. Oddness determines the variation
at the opposite centre. If $B\leq p$, this gives $B-1=E$
parameters. If $B>p$, it gives $(p-2)+(b_2-1)=E-1$ parameters.
For the last parameter vary $\xi=\sqrt b$. If $h$ is the coefficient
of $(t-\xi)^{-(p-2)}$, differentiation with respect to $\xi$
produces the term
\[
 (p-2)h(t-\xi)^{-(p-1)}.
\]
No change of a coefficient produces a pole of this order. The
change of centre is therefore independent of the other changes,
which are themselves independent by uniqueness of partial fractions.

To compute their effect at infinity, set
\[
 D_0(z)=
 \begin{cases}
 (z-b)^E & B\leq p,\\
 (z-b)^{p-1}(z-c_2)^{b_2-1} & B>p .
 \end{cases}
\]
In both cases $\deg D_0=E$. The $E$ changes of $P_{\mathrm{aux}}$
are odd, vanish at infinity, and have denominator dividing
$D_0(t^2)$. They therefore belong to the space
\begin{equation}\label{eq-directions}
 \left\{\frac{tL(t^2)}{D_0(t^2)}
        \ \middle|\ L\in k[z],\ \deg L<E\right\}.
\end{equation}
This space has dimension $E$. Since the changes are independent,
they form a basis. We need only show that an element of
\eqref{eq-directions} whose coefficients in
\eqref{eq-corr-aux-rows} vanish is zero.

Set $s=t^{-2}$ and $D_*(s)=s^ED_0(1/s)=\sum_{i=0}^E d_i s^i$.
An element of \eqref{eq-directions} has the form
\[
 t^{-1}\frac{A(s)}{D_*(s)}
   =t^{-1}\sum_{j\geq0}u_js^j,
 \deg A<E.
\]
Thus $u_j$ is the coefficient of $t^{-(2j+1)}$. Multiplication by
$D_*$ gives
\[
 \sum_{i=0}^E d_i u_{j+E-i}=0
 \quad\text{for every }j\geq0.
\]
Both $d_0=1$ and $d_E$ are nonzero. Hence $E$ consecutive values
of $u_j$ determine the sequence both forwards and backwards.

Set $j_0=(K-1)/2$. Before we omit $pa$, the coefficients in
\eqref{eq-corr-aux-rows} are the $E+1$ consecutive terms
$u_{j_0},\ldots,u_{j_0+E}$. Since $pa=K+2a$, the omitted term
is $u_{j_0+a}$. In the recurrence at $j=j_0$ its coefficient
is $d_{E-a}=d_{e-1}$. If $B\leq p$, then
\[
 d_{e-1}=\binom E{e-1}(-b)^{e-1}\ne0
\]
because $0\leq e-1\leq E<p$. If $B>p$, then
$D_*(s)=(1-bs)^{p-1}(1-c_2s)^{b_2-1}$. Up to degree $e-1<p$,
the first factor has the same coefficients as $(1-bs)^{-1}$.
Consequently,
\begin{equation}\label{eq-recurrence}
 d_{e-1}
 =[s^{e-1}]\frac{(1-c_2s)^{b_2-1}}{1-bs}
 =b^{e-b_2}(b-c_2)^{b_2-1}\ne0 ,
\end{equation}
where $e-b_2=p-a-1\geq0$.

If the retained coefficients vanish, the recurrence forces the
omitted coefficient to vanish as well. It then gives $u_j=0$ for
every $j$, so $A=0$. Thus the map from the $E$ parameter changes
to the $E$ coefficient changes has zero kernel. Its Jacobian is
invertible.
\end{proof}

We now give the two factors of $G$. First suppose $e>0$. Choose
one centre from each pair of nonzero poles of $P$, and let $A(t)$
be the sum of the principal parts at the chosen centres. Since $P$
is odd, $P_{\mathrm{aux}}(t)=A(t)-A(-t)$. Lift the centres and
coefficients of $A$ to $R$, denote the resulting function by
$\widehat A$, and set
\begin{equation}\label{eq-auxfactor}
 G_h(T)=
 \frac{1+\varpi_{s_h}\widehat A(T/\varpi_{r_h})}
      {1+\varpi_{s_h}\widehat A(-T/\varpi_{r_h})}.
\end{equation}
In the variable $T/\varpi_{r_h}$, the zeros and poles of the
numerator reduce to the chosen centres. Those of the denominator
reduce to their negatives. These two sets are disjoint, so there is
no cancellation. By Lemma~\ref{lem-auxiliary-branches}, this gives
$2B$ branch points, all of valuation $r_h$. Moreover, on writing
$T=\varpi_{r_h}U$ and $t=\bar U$, we have
\[
 \overline{\frac{G_h(\varpi_{r_h}U)-1}{\varpi_{s_h}}}
 =A(t)-A(-t)=P_{\mathrm{aux}}(t).
\]
This last equality is why we used the principal parts of $P$ to
define $G_h$.

For the points of valuation $r_c$, take the polynomial $f$ of
Lemma~\ref{lem-data}, scaled to the coefficient of the fixed
quotient. Use the polynomial of degree $N=(p-3)a$ if $e>0$, and
the one of degree $N=(p-1)a$ if $e=0$. Set $n=N/2$. From each
pair of roots choose one, denoted $\bar x_i$, and choose a lift
$x_i\in R^\times$. Let $R_i\in\mathbb Z$ represent
\[
 \operatorname{Res}_{t=\bar x_i}\frac{dt}{t^{a+1}f(t)}
 =\frac{1}{\bar x_i^{a+1}f'(\bar x_i)}
 \in\mathbb F_p^\times.
\]
Set
\begin{equation}\label{eq-critfactor}
 G_c(T)=\prod_{i=1}^n
 \left(\frac{T-\varpi_{r_c}x_i}{T+\varpi_{r_c}x_i}\right)^{R_i}.
\end{equation}
This gives $N$ distinct points $T=\pm\varpi_{r_c}x_i$. When
$N=0$, the product is $1$.

Finally, take $G=G_cG_h$ if $e>0$ and $G=G_c$ if $e=0$.
Each factor satisfies $G_\bullet(-T)=G_\bullet(T)^{-1}$. At every
new point the order of $G$ is prime to $p$, so its ramification
index in $Y^{p^2}=A_0G^p$ is $p$. At each of the $a+1$ old points
the order of $A_0$ is prime to $p$, so the index is $p^2$. Thus
the fixed quotient, the reflection, and all branch counts are as
required. It remains to choose the coefficients of $G$ so that the
reduction is smooth and its second Witt coordinate is $f_2$. This
is the purpose of the next two sections.

\section{Adjusting the branch points}\label{sec-corr-formal}

In the preceding section we chose the branch points of the cover
\[
 Y^{p^2}=A_0(T)G(T)^p.
\]
We now vary these points without changing their reductions. The
purpose is to prescribe a term in the equation at a given radius
$r$. When both factors $G_c$ and $G_h$ are present, we must also
preserve the term already prescribed at $r_h$. These are the two
requirements in \eqref{eq-correction} below.

We first explain which quantities may be varied. In $G_c$ we
replace each $x_i$ by $x_i+z_i$. In $G_h$ we vary the $E$
coefficients and centres of Lemma~\ref{lem-auxiliary}, and denote
their increments by $y_1,\ldots,y_E$. We require all increments to
belong to the maximal ideal of the coefficient ring. Thus $x_i+z_i$
and $x_i$ have the same reduction, and the same holds for each
centre and coefficient in $G_h$. Denote the resulting factors by
$\widetilde G_c$ and $\widetilde G_h$.

There is one further choice that does not move any branch point.
Replacing $G$ by $\widetilde G$ changes the character of the cover
by the degree-$p$ character $\kappa_1(\widetilde G/G)$. Indeed,
\[
 \kappa_2(A_0\widetilde G^p)-\kappa_2(A_0G^p)
 =\kappa_1(\widetilde G/G),
\]
where we identify $\mathbb Z/p\mathbb Z$ with the subgroup of order
$p$ in $\mathbb Z/p^2\mathbb Z$. We may therefore multiply
$\widetilde G/G$ by $H^p$ without changing this character. We use
the branch points to prescribe certain coefficients and the factor
$H^p$ to prescribe the remaining ones.

Recall that $v(p)=1$, that $v(\varpi_s)=s$, and that the elements
$\varpi_s$ are chosen compatibly under multiplication. We write
$v_r$ for the Gauss valuation in the variable $T/\varpi_r$. Thus
the term $bT^{-q}$ has valuation $v(b)-qr$ at $r$.

\begin{lemma}\label{lem-correction}
Suppose $e>0$ and let $0\leq r<r_h$ be rational. Let $F$ be an
odd polynomial in $T^{-1}$ such that $\deg_{T^{-1}}F<M$ and
$v_r(F)=0$. Set $\beta=\alpha-Mr$, where $\alpha=p/(p-1)$.
After a finite extension of the coefficient field, one can choose
the increments above and a rational function $H$ satisfying
$H(-T)=H(T)^{-1}$ so that
\[
 Q=\frac{\widetilde G_c}{G_c}
   \frac{\widetilde G_h}{G_h}H^p
\]
satisfies
\begin{equation}\label{eq-correction}
 v_r\bigl(Q-(1-\varpi_\beta F)\bigr)>\beta
 \quad\text{and}\quad
 v_{r_h}(Q-1)>s_h.
\end{equation}

There is also an assertion using only $G_c$. Let $N$ be the degree
of either polynomial in Lemma~\ref{lem-data} and set $K=N+a$. Let
$0\leq r<r_c$ be rational and let $F$ be an odd polynomial in
$T^{-1}$ with $\deg_{T^{-1}}F<K$ and $v_r(F)=0$. Set
$\beta=K(r_c-r)$. Then one can choose the root increments and $H$
so that
\[
 v_r\bigl((\widetilde G_c/G_c)H^p-(1-\varpi_\beta F)\bigr)>\beta .
\]
For $F=0$, both assertions hold with all increments zero and $H=1$.
\end{lemma}

\begin{proof}
We first prove the assertion with both factors. Write
$n=N/2$, where $N=(p-3)a$, and recall that
\[
 K=N+a,\quad
 E=a+e-1,\quad
 \Delta=r_c-r_h,\quad
 s_h=K\Delta.
\]
For a rational function $U$ regular at infinity, write
$c_q(U)=[T^{-q}]U$ for the coefficient of $T^{-q}$ in its
expansion there. We shall first match the coefficients in degrees
$q<M$. We shall then bound the terms in degrees $q\geq M$.

The symmetry $U(-T)U(T)=1$ reduces the number of coefficients to
be prescribed by half. Indeed, if $U=1+\sum_{q\geq1}b_qT^{-q}$,
then for every even $q>0$ we have
\[
 2b_q=-\sum_{i=1}^{q-1}(-1)^i b_i b_{q-i}.
\]
Since $p$ is odd, the even coefficients are determined successively
by the odd ones. We therefore prescribe only the coefficients in
odd degrees.

We next determine the power of $\varpi$ by which to divide each
coefficient before reducing it. In degree $q>0$, the factor
$\widetilde G_c/G_c$ contributes a multiple of $\varpi_{qr_c}$.
The factor $\widetilde G_h/G_h$ contributes a multiple of
$\varpi_{s_h+qr_h}$ by \eqref{eq-auxfactor}. Set
\begin{equation}\label{eq-corr-weight}
 w(q)=\min\{qr_c,s_h+qr_h\}
     =qr_h+\Delta\min\{q,K\}.
\end{equation}
Thus $w(q)=qr_c$ for $q\leq K$ and $w(q)=s_h+qr_h$ for $q\geq K$.
In particular, $w(M)=s_h+Mr_h=\alpha$. We shall use the inequality
\begin{equation}\label{eq-corr-subadd}
 w(i)+w(j)\geq w(i+j),
\end{equation}
which follows from
$\min\{i,K\}+\min\{j,K\}\geq\min\{i+j,K\}$.

For each odd $\ell$ with $1\leq\ell\leq a$, introduce a parameter
$h_\ell$. Set
\begin{equation}\label{eq-Hh}
 Z=\sum_{\substack{1\leq\ell\leq a\\ \ell\ {\rm odd}}}
       h_\ell\varpi_{w(p\ell)/p}T^{-\ell}
 \quad\text{and}\quad
 H=\frac{1+Z}{1-Z}.
\end{equation}
Then $H(-T)=H(T)^{-1}$. The reason for introducing $h_\ell$ is that
$H^p$ contributes $2h_\ell^p$ to the reduced equation in degree
$p\ell$, once the parameters with smaller index have been set to
zero. This will allow us to treat the degrees divisible by $p$.

For each odd $q<M$, set $\Phi_q=c_q(Q-1)/\varpi_{w(q)}$. There are
as many parameters as equations, since
\begin{equation}\label{eq-corr-count}
 n+E+\frac{a+1}{2}=\frac{M-1}{2}.
\end{equation}
Let $R$ denote the coefficient ring, enlarged if necessary to contain
the finitely many elements $\varpi_s$ used here. We claim that
$\Phi_q$ is a power series over $R$ in the parameters $z_i,y_j,h_\ell$.
Only the factor $H^p$ needs to be checked. Consider a monomial
$\prod_\ell h_\ell^{m_\ell}$ in degree $q=\sum_\ell\ell m_\ell$.
If every $m_\ell$ is divisible by $p$, its coefficient has valuation
at least $\sum_\ell(m_\ell/p)w(p\ell)\geq w(q)$ by
\eqref{eq-corr-subadd}. Otherwise its coefficient before inserting
the powers of $\varpi$ in \eqref{eq-Hh} is divisible by $p$.
Indeed, modulo $p$, a $p$-th power contains only monomials whose
exponents are all divisible by $p$. In this case the valuation is
at least
\begin{equation}\label{eq-corr-mixed}
 1+\frac1p\sum_\ell m_\ell w(p\ell)
 \geq 1+\frac{w(pq)}p
 \geq 1+\frac{w(q)}p
 >w(q),
\end{equation}
where the last inequality follows from $w(q)<\alpha$. Products
with terms from the other factors preserve these bounds by
\eqref{eq-corr-subadd}. This proves the claim. Also, $\Phi_q(0)=0$
since $Q=1$ when all parameters are zero.

We now show that the $\Phi_q$ can be prescribed. By
Lemma~\ref{lem-formal}, it suffices to prove that
\[
 \mathcal B=
 k[[z_1,\ldots,z_n,y_1,\ldots,y_E,(h_\ell)_\ell]]
 \big/
 \bigl(\overline{\Phi}_q\mid 1\leq q<M,\ q\ {\rm odd}\bigr)
\]
has finite positive dimension over $k$. We prove this by showing
that every parameter is nilpotent in $\mathcal B$. Work in
$\mathcal B_{\mathrm{red}}$, the quotient obtained by setting all
nilpotent elements equal to zero.

We first determine the $z_i$. Set $X=\varpi_{r_c}/T$. In degrees
$q<K$, division by $\varpi_{w(q)}=\varpi_{qr_c}$ amounts to
taking the coefficient of $X^q$. In these degrees the reduction of
$\widetilde G_h/G_h$ is $1$, since
\[
 s_h+qr_h-qr_c=(K-q)\Delta>0.
\]
The odd coefficients of $Q-1$ vanish by the defining equations of
$\mathcal B$. The recursion for the even coefficients therefore
gives $Q\equiv1\pmod{X^K}$ after reduction.

Let $W$ be the reduction of $\widetilde G_c/G_c$ in the variable
$X$. By \eqref{eq-corr-mixed}, the reduction of $H^p$ through
degree $K-1$ contains only powers of $X^p$. Its derivative is
therefore zero in these degrees when taken with $X\,d/dX$.
Logarithmically differentiating the preceding congruence gives
$XW'/W\equiv0\pmod{X^K}$. For a single root $x$, the calculation is
\[
 \frac{X}{(1-xX)/(1+xX)}
 \frac{d}{dX}\left(\frac{1-xX}{1+xX}\right)
 =-\frac{2xX}{1-x^2X^2}
 =-2\sum_{\substack{q\geq1\\q\ {\rm odd}}}x^qX^q.
\]
Write $\bar x_i$ for the original reduced roots and
$r_i\in\mathbb F_p^\times$ for their residues. Comparing the
coefficients for the old and new roots gives
\[
 \sum_{i=1}^n r_i(\bar x_i+z_i)^q
 =\sum_{i=1}^n r_i\bar x_i^q
 \quad\text{for odd }q<K\text{ with }p\nmid q.
\]
These are precisely the moment equations \eqref{eq-moments} at the
original tuple $(\bar x_i)$. By Lemma~\ref{lem-finitecartier},
their completed local ring at this tuple is finite dimensional
over $k$. Its maximal ideal is therefore nilpotent. The elements
$z_i$ are images of this ideal in $\mathcal B_{\mathrm{red}}$,
which has no nonzero nilpotent elements. Hence $z_i=0$ for every
$i$.

We next determine $h_\ell$ for odd $\ell\leq a-2$. For these
indices $p\ell<K$, so $w(p\ell)=p\ell r_c$. The critical ratio is
now $1$ since all $z_i$ vanish, and the auxiliary ratio is $1$
through degree $K-1$ after reduction. In degree $p\ell$, the
reduced coefficient of $H^p-1$ is $2h_\ell^p$ plus terms involving
some $h_j$ with $j<\ell$. This follows by expanding
$(1+Z)^p/(1-Z)^p$ and using \eqref{eq-corr-mixed}. Taking the
indices $\ell$ in increasing order gives $2h_\ell^p=0$ and hence
$h_\ell=0$ in $\mathcal B_{\mathrm{red}}$.

It remains to determine $y_1,\ldots,y_E$ and $h_a$. For
$K\leq q<M$, the scale is $w(q)=s_h+qr_h$. In the auxiliary ratio,
every term containing two or more powers of $\varpi_{s_h}$ vanishes
after division by $\varpi_{w(q)}$ and reduction. Thus the
contribution of this ratio is
\[
 [t^{-q}](\widetilde P_{\mathrm{aux}}-P_{\mathrm{aux}}),
\]
where $\widetilde P_{\mathrm{aux}}$ is obtained by varying the
centres and coefficients by $y_1,\ldots,y_E$. The only contribution
from $H^p-1$ is $2h_a^p$ in degree $pa$. Indeed, a term in degree
$jpa$ with $j\geq2$ has valuation at least
$jw(pa)=js_h+jpa r_h>w(jpa)$. A product of the degree-$pa$ term
with a nonconstant auxiliary term also has an extra factor
$\varpi_{s_h}$ and vanishes after reduction.

Take first the odd degrees $K\leq q<M$ with $q\ne pa$. These give
$E$ equations involving only $y=(y_1,\ldots,y_E)$. Denote their
left hand sides by $\Psi(y)$. Since $\Psi(0)=0$, we may write
$\Psi(y)=A(y)y$, where $A(0)$ is the Jacobian in
Lemma~\ref{lem-auxiliary}. Its determinant is nonzero. Thus
$\det A(y)$ is a unit in the local ring, and $\Psi(y)=0$ implies
$y=0$. The equation in degree $pa$ now reads $2h_a^p=0$, whence
$h_a=0$.

We have shown that every parameter vanishes in
$\mathcal B_{\mathrm{red}}$. Thus every parameter is nilpotent in
$\mathcal B$. There are finitely many parameters, so their
monomials span a finite dimensional $k$-space. Since the origin
satisfies all the equations, $\mathcal B\ne0$. Lemma~\ref{lem-formal}
now shows that the $\Phi_q$ may be assigned any values of positive
valuation, after a finite extension, with all parameter increments
of positive valuation.

We apply this conclusion to the prescribed polynomial $F$. The
function $1-\varpi_\beta F$ need not satisfy the reflection
identity. We therefore use
\begin{equation}\label{eq-corr-target}
 u=\varpi_\beta F
 \quad\text{and}\quad
 J=\frac{1-u/2}{1+u/2}.
\end{equation}
Since $F$ is odd, $J(-T)J(T)=1$. Moreover,
$J-(1-u)=u^2/(2+u)$, so $v_r(J-(1-u))=2\beta>\beta$.
It is consequently enough to match $J$ through degree $M-1$ and
to bound the remaining terms of $Q-1$.

For $q>0$, the coefficient $c_q(J-1)$ has valuation at least
$\beta+qr$. Indeed, $v_r(u)=\beta$, and each higher power $u^j$
has valuation $j\beta$ at $r$. For $q<M$ we have
\begin{equation}\label{eq-excess}
 \beta+qr-w(q)
 =(M-q)(r_h-r)+\max\{K-q,0\}\Delta>0.
\end{equation}
We may therefore prescribe
$\Phi_q=c_q(J-1)/\varpi_{w(q)}$ for every odd $q<M$. The chosen
$Q$ and $J$ have constant term $1$ and both satisfy the reflection
identity. The recursion for the even coefficients gives
$c_q(Q)=c_q(J)$ for every $q<M$.

We finish by bounding the coefficients for $q\geq M$. Set
\begin{equation}\label{eq-corr-mu}
 \mu=\frac{s_h}{pa}.
\end{equation}
Since $K<pa<M$ and $s_h=K\Delta$, we have $M\mu>s_h$ and
$\mu<\Delta$. For every $\ell\leq a$,
\begin{equation}\label{eq-corr-slopes}
 \frac{w(p\ell)}p
 =\min\{\ell r_c,s_h/p+\ell r_h\}
 \geq \ell(r_h+\mu).
\end{equation}
Choose $\epsilon>0$ no larger than the valuation of any nonzero
increment, including the $h_\ell$. Set
$S=(\widetilde G_c/G_c)H^p$ and $U=\widetilde G_h/G_h$, so
$Q=SU$. In degree $q>0$, the coefficients satisfy
\begin{equation}\label{eq-tail}
 \begin{aligned}
 v\bigl(c_q(S-1)\bigr)
     &\geq q(r_h+\mu)+\epsilon,\\
 v\bigl(c_q(S(U-1))\bigr)
     &\geq s_h+qr_h+\epsilon.
 \end{aligned}
\end{equation}
For the first bound, use \eqref{eq-corr-slopes} for $H$ and
$r_c>r_h+\mu$ for the root factor. Every nonconstant term contains
at least one increment. For the second bound, the coefficient of
degree $j>0$ in $U-1$ already contains $\varpi_{s_h+jr_h}$ and
an auxiliary increment. Multiplying by a term of degree $q-j$
in $S$ adds valuation at least $(q-j)r_h$.

Since $Q-1=(S-1)+S(U-1)$ and $M\mu>s_h$, it follows for $q\geq M$
that
\[
 v\bigl(c_q(Q-1)\bigr)-qr
 \geq \beta+\epsilon+(q-M)(r_h-r)>\beta .
\]
This bound tends to infinity with $q$. In degrees below $M$, the
coefficients of $Q$ agree with those of $J$, and
$v_r(J-(1-u))=2\beta$. In degrees at least $M$, the polynomial
$1-u$ has no terms. This proves the first inequality of
\eqref{eq-correction}.

For the estimate at $r_h$, set $V=\varpi_{r_h}/T$. If $0<q<M$,
the coefficient of $V^q$ in $Q-1$ agrees with that in $J-1$ and
has valuation at least
\[
 \beta-q(r_h-r)=s_h+(M-q)(r_h-r)>s_h.
\]
For $q\geq M$, the same coefficient has valuation at least
$s_h+\epsilon$ by \eqref{eq-tail}. Thus every coefficient of the
formal expansion in $V$ has valuation greater than $s_h$.

The auxiliary factor may have poles on $v(T)=r_h$, so we pass
from the formal expansion to the rational function by clearing
denominators. Write $Q-1=P(V)/D(V)$ with $P,D\in R[V]$ and
$D(0)=1$. Such a denominator exists by the formulas for $G_c$,
$G_h$ and $H$. In particular, an auxiliary term
$(T/\varpi_{r_h}-\xi)^{-j}$ becomes $V^j/(1-\xi V)^j$.
Multiplication of the formal expansion by $D$ shows that every
coefficient of $P$ has valuation greater than $s_h$. Since $D$
has Gauss valuation zero, $v_{r_h}(Q-1)>s_h$. This proves the
second inequality of \eqref{eq-correction}.

Finally, suppose that only $G_c$ is used. Take $w(q)=qr_c$ and
retain only the parameters $z_i$ and $h_\ell$ with $\ell$ odd and
$\ell\leq a-2$. For either polynomial of Lemma~\ref{lem-data},
$p(a-2)<K\leq pa$, so these are precisely the odd positive
indices for which $p\ell<K$. The number of parameters is
$n+(a-1)/2=(K-1)/2$. In the variable $X=\varpi_{r_c}/T$,
both $\widetilde G_c/G_c$ and $H$ have integral coefficients.
The same logarithmic derivative gives the moment equations, so
the $z_i$ vanish in the reduced ring. The coefficients in degrees
$p\ell<K$ then give $h_\ell=0$ in increasing order of $\ell$.
Lemma~\ref{lem-formal} therefore applies to the odd coefficients
below $K$.

Use $J$ from \eqref{eq-corr-target} with $\beta=K(r_c-r)$.
For $q<K$ the target coefficient has excess valuation at least
$\beta+qr-qr_c=(K-q)(r_c-r)>0$. Thus $Q$ and $J$ can be made to
agree in every degree below $K$. For $q\geq K$ we have
$v(c_q(Q-1))\geq qr_c+\epsilon$, since each nonconstant term
contains an increment. After subtracting $qr$, this gives at least
$K(r_c-r)+\epsilon>\beta$. The first inequality follows as before.
\end{proof}

We can therefore prescribe the term $-\varpi_\beta F$ at $r$
while keeping the reductions of the branch points fixed. When both
$G_c$ and $G_h$ are used, the inequality $v_{r_h}(Q-1)>s_h$ also
keeps the previously chosen term of valuation $s_h$ at $r_h$
unchanged. Section~\ref{sec-boundedproof} uses the first property
to correct the equation and the second to preserve the correction
already made at $r_h$.
\section{Lifting actions with bounded jumps}\label{sec-boundedproof}

Let $p$ be an odd prime. Let $1\leq a<p$ be odd, let $0\leq e<p$,
and set $M=pa+2e$. We prove that a local $D_{p^2}$-action with
cyclic upper jumps $(a,M)$ lifts to characteristic zero. By
Lemma~\ref{pre-normalized-witt}, we may write its Witt vector as
$(x^{-a},f_2)$, where $f_2$ is an odd standard polynomial in
$x^{-1}$. We fix the lift $\chi_1$ of the first equation
$w^p-w=x^{-a}$ given by \eqref{eq-firstdepth}.

There are two things to prove. First, we must construct a cover
with good reduction and jumps $(a,M)$. Second, we must show that
its reduction has the prescribed second coordinate $f_2$. These
are distinct requirements. For instance, when $e>0$, the jump $M$
tells us that $\deg_{x^{-1}}f_2=M$, but it does not determine the
coefficient of $x^{-M}$.

We use the notation of the preceding two sections. In particular,
$\alpha=p/(p-1)$, $r_c=1/((p-1)a)$, and $K=N+a$. The choices in
the two cases are
\[
 \begin{array}{c|c|c|c}
       & N & K & r_* \\ \hline
 e=0   & (p-1)a & pa     & 0\\
 e>0   & (p-3)a & (p-2)a & r_h .
 \end{array}
\]
Here $N$ is the number of branch points supplied by $G_c$. If
$e>0$, the factor $G_h$ supplies another $2B=M-K$ branch points,
where $B=a+e$. For a character $\psi$, write $\delta_\psi(r)$ for
its depth on $D[r]=\{v(T)\geq r\}$. When the depth is positive,
write $\omega_\psi(r)$ for its differential in the variable
$t=\overline{T/\varpi_r}$.

We use the following elementary observation twice. Suppose that a
continuous piecewise linear function $\delta$ has slopes at most
$d$ on $[u,v]$. If $\delta(v)-\delta(u)=d(v-u)$, then every slope
on this interval equals $d$. Thus, after correcting the differential
at $u$, it suffices to preserve the depths at $u$ and $v$. No
separate calculation at the intermediate radii is needed.

\begin{lemma}\label{lem-bounded-critical-depth}
Choose the polynomial of degree $N$ in Lemma~\ref{lem-data}, scaled
to the coefficient of the fixed quotient. After a finite extension
of the coefficient field, its roots can be lifted in $G_c$ so that
$\chi_c=\kappa_2(A_0G_c^p)$ satisfies
\[
 \delta_{\chi_c}(r)=\alpha-K(r_c-r)
 \quad\text{for }r_*\leq r\leq r_c.
\]
If $e>0$, the coefficient of $t^{-K-1}dt$ in
$\omega_{\chi_c}(r_h)$ is nonzero. It is the same for all such
choices of $G_c$ with the prescribed ordered root reductions.
\end{lemma}

\begin{proof}
We first compute the differential at $r_c$. Set
$\chi_0=\kappa_2(A_0)$, so that $p\chi_0=\chi_1$. By
\eqref{eq-firstdepth}, the quotient has depth
$ar_c=1/(p-1)$ and differential $c_*t^{-a-1}dt$ for some
$c_*\in k^\times$. It follows from
\cite[Theorem~4.3(iii)]{Wewers} that $\delta_{\chi_0}(r_c)=\alpha$
and
\[
 \Cart(\omega_{\chi_0})=\omega_{\chi_0}+c_*t^{-a-1}dt
\]
at $r_c$. Since $c_*t^{-a-1}dt=d(-(c_*/a)t^{-a})$, the right
hand side is fixed by Cartier. It therefore has only simple poles.
All branch points of $A_0$ have valuation greater than $r_c$, so
there is no pole away from $0$ and $\infty$, by
\cite[Proposition~6.19]{Obus}. A differential on $\mathbb P^1_k$
with at most simple poles at $0$ and $\infty$ is a multiple of
$dt/t$. The reflection $t\mapsto-t$ fixes $dt/t$ and negates the
differential above. Since $p$ is odd, that multiple is zero.

Let $g$ be the reduction of $G_c$ in the variable $T/\varpi_{r_c}$.
Adding the character of $G_c$ gives
\begin{equation}\label{eq-criticalbase}
 \begin{aligned}
 \omega_{\chi_0}(r_c)&=-c_*t^{-a-1}dt,\\
 \omega_{\chi_c}(r_c)&=-c_*t^{-a-1}dt+\frac{dg}{g}.
 \end{aligned}
\end{equation}
The term $dg/g$ is regular at zero, so the second differential is
nonzero and $\delta_{\chi_c}(r_c)=\alpha$. For the prescribed
roots, Lemma~\ref{lem-data} identifies it with
$dt/(t^{a+1}f(t))$, where $\deg f=N$. Its order at infinity is
$a+N-1=K-1$. Hence the left slope at $r_c$ is $K$ by
\eqref{eq-slopes}.

We now choose lifts of the roots for which this slope persists down
to $r_*$. We give the parameter space explicitly because
Theorem~\ref{thm-minimum} requires a quasi-compact family. Write
one point from each opposite pair as $\varpi_{r_c}y_i$, where
$|y_i|=1$, and require $|y_i^2-y_j^2|=1$ for $i\ne j$.
By Lemma~\ref{lem-finitecartier}, the moment equations
\eqref{eq-moments} have only finitely many ordered solutions in
the reduction of this space. Choose a polynomial over $k$ that is
$1$ at the prescribed solution and $0$ at all the others. Lift it
to a polynomial $P$ over the valuation ring and impose $|P(y)|=1$.
The resulting space $\mathcal A_c$ is affinoid. It contains every
lift of the prescribed tuple, and no other solution of the moment
equations occurs as a reduction in $\mathcal A_c$. If $N=0$, take
$\mathcal A_c$ to be a point.

We verify the remaining hypotheses of Theorem~\ref{thm-minimum}.
The $a+1$ branch points of $A_0$ have index $p^2$, and the $N$
points of $G_c$ have index $p$. They are disjoint throughout
$\mathcal A_c$. Locally along a branch point of index $d$, the
normalized map has equation $u\mapsto u^d$ after an \'etale base
change, since the unit in the Kummer equation has a $d$-th root
after such a change. Thus both maps in the tower are finite flat
maps of smooth relative curves, and their branch points form
finite \'etale families.

Choose $s_1>r_c$ smaller than the valuation of every old branch
point. For $0<r<s_1$, the inverse image of $D(r)=\{v(T)>r\}$ is
connected. Indeed, the group acts transitively on its components,
and a point with full inertia lies on a component fixed by the
whole group. For rational $0<r\leq r_c$, the quotient has positive
depth $ar$. After the extension of constants allowed in
Section~\ref{sec-analytic}, its residue extension has purely
inseparable degree $p$. The decomposition group of the full cover
maps onto the group of the quotient, so it is $C_{p^2}$. If the
residue extension had a nontrivial separable part, the henselian
extension would have an unramified subextension of degree $p$.
Its unique degree-$p$ subextension is the fixed quotient, which
has positive depth. This is a contradiction. Finally, the inverse
image of $D[r]$ in the quotient is a closed disk by
\eqref{eq-firstdepth} and \cite[Corollary~2.18]{Kinks}. All $K+1$
branch points lie in $D[r_c]$. The theorem therefore applies with
$r_0=r_c$ and gives the slope bound $K$.

For each choice of roots, set
\[
 \lambda=\sup\Bigl(\{r_*\}\cup
 \{r\in(r_*,r_c]\cap\mathbb Q
       \mid \delta_{\chi_c}'(r^-)<K\}\Bigr).
\]
Thus $\lambda$ is the lower endpoint of the interval ending at
$r_c$ on which the depth has the required slope. By
Theorem~\ref{thm-minimum}, we may choose a fibre minimizing
$\lambda$, after a finite extension. The prescribed reductions
give slope $K$ on an interval ending at $r_c$, so the minimum is
less than $r_c$. Conversely, by \eqref{eq-criticalbase} and
\eqref{eq-log}, slope $K$ at $r_c$ implies the moment equations.
The minimizing fibre therefore has the prescribed root reductions.

Suppose $\lambda>r_*$. Then $\lambda$ is a rational breakpoint.
On $[\lambda,r_c]$ we have
\begin{equation}\label{eq-criticalline}
 \begin{aligned}
 d_c(r):=\delta_{\chi_c}(r)&=\alpha-K(r_c-r),\\
 d_c(r)-par&=(pa-K)(r_c-r).
 \end{aligned}
\end{equation}
All branch points have valuation greater than $\lambda$. Thus
$\omega_{\chi_c}(\lambda)$ is regular away from zero, including at
infinity, as in Section~\ref{sec-analytic}. Its pole at zero has
order $K+1$, since its right slope is $K$. We claim that
\[
 \omega_{\chi_c}(\lambda)=C_\lambda t^{-K-1}dt+df,
\]
where $C_\lambda\ne0$ and $f$ is an odd polynomial in $t^{-1}$
of degree less than $K$.

If $e>0$, then $K=(p-2)a$ and $d_c(\lambda)>pa\lambda$.
Equation~\eqref{eq-cartierdepth} shows that the differential is
exact. Its leading term is also exact since $p\nmid K$.
If $e=0$, then $K=pa$ and $d_c(\lambda)=p\delta_{\chi_1}(\lambda)$.
In this case
$\Cart(\omega_{\chi_c}(\lambda))=\omega_{\chi_1}(\lambda)$.
Both slopes of the quotient equal $a$, and its differential has
no pole away from zero. Thus
$\omega_{\chi_1}(\lambda)=ct^{-a-1}dt$ for some $c\ne0$.
Since
\[
 \Cart(C_\lambda t^{-pa-1}dt)
   =C_\lambda^{1/p}t^{-a-1}dt,
\]
we obtain $C_\lambda=c^p$. The difference again has Cartier image
zero and is exact. In either case, integrating the remaining terms
gives a polynomial of degree less than $K$. Replacing it by
$(f(t)-f(-t))/2$ preserves its derivative and makes it odd.

If $df\ne0$, write $f=\sum_j f_jt^{-j}$, lift its nonzero
coefficients to units $\widetilde f_j$, and set
\[
 F=\sum_j\widetilde f_j(\varpi_\lambda/T)^j
 \quad\text{and}\quad
 \beta=K(r_c-\lambda).
\]
Then $v_\lambda(F)=0$ and $F$ reduces to $f$. The single-factor
assertion of Lemma~\ref{lem-correction} gives a correction $Q$ with
$v_\lambda(Q-(1-\varpi_\beta F))>\beta$. By
\eqref{eq-kummerdepth}, its character has depth
$\alpha-\beta=d_c(\lambda)$ and differential $-df$. Thus the new
differential at $\lambda$ is $C_\lambda t^{-K-1}dt$. If $df=0$,
no change is needed.

The depth at $\lambda$ is unchanged because the remaining
differential is nonzero. The depth at $r_c$ is also unchanged.
Their difference is $K(r_c-\lambda)$, so the slope remains $K$
throughout $[\lambda,r_c]$. Moreover,
$C_\lambda t^{-K-1}dt$ has order $K-1$ at infinity. The left
slope at $\lambda$ is therefore $K$, and the interval of slope
$K$ extends below $\lambda$. This contradicts minimality. Hence
$\lambda=r_*$, and continuity gives the asserted depth at $r_*$.

It remains to prove independence of the coefficient when $e>0$.
Take two choices $G_{c,1},G_{c,2}$ satisfying the first assertion
with the same root reductions. Set $\Delta=r_c-r_h$ and
$U=T/\varpi_{r_h}$. Since the ratios of the root factors reduce
to $1$ in the variable $\varpi_{r_c}/T$, we have
\begin{equation}\label{eq-criticalraw}
 G_{c,2}/G_{c,1}=1+\sum_{i>0}b_iU^{-i},
 \quad\text{with }v(b_i)>i\Delta .
\end{equation}
In particular, $v(b_i)>s_h=K\Delta$ for $i\geq K$. The character
of this ratio has depth at most $Mr_h$ at $r_h$. If its depth is
smaller, the two differentials agree. Otherwise, apply
Lemma~\ref{an-cutoff} with $m=K$ and $\beta=s_h$, using the
normalized representative from Section~\ref{sec-analytic} as
$F_1$ and the ratio \eqref{eq-criticalraw} as $F_2$. This is
permitted since $p\nmid K$ and $\alpha-s_h=Mr_h$. The normalized
representative also has coefficients of valuation greater than
$s_h$ in degrees $i\geq K$. Its reduced differential therefore
has zero coefficient of $t^{-K-1}dt$. The two original
coefficients are equal. Their common value is nonzero, since the
right slope at $r_h$ is $K$.
\end{proof}

When $e=0$, the preceding lemma already gives
$\delta_{\chi_c}(r)=par$ down to $r=0$. For $e>0$, it gives the
required depth only as far as $r_h$. We now use the $2B$ branch
points of $G_h$ to obtain slope $M$ from $r_h$ down to zero.

\begin{lemma}\label{lem-bounded-both-depths}
Suppose $e>0$. Let $C$ be the coefficient determined by
Lemma~\ref{lem-bounded-critical-depth}, and let $\eta$ be the
differential of Lemma~\ref{lem-auxiliary} with leading term
$Ct^{-K-1}dt$ at zero. After a finite extension, one can choose
$G_c$ and $G_h$ so that $\chi=\kappa_2(A_0(G_cG_h)^p)$ satisfies
\[
 \omega_\chi(r_h)=\eta
 \quad\text{and}\quad
 \delta_\chi(r)=Mr\ \text{ for }0\leq r\leq r_h.
\]
The resulting cover has good reduction with jumps $(a,M)$.
\end{lemma}

\begin{proof}
Choose $G_c$ by Lemma~\ref{lem-bounded-critical-depth}, and let
$\eta_0$ denote the principal part of $\eta$ at zero. Since
$\delta_{\chi_c}(r_h)=Mr_h>par_h$, the differential
$\omega_{\chi_c}(r_h)$ is exact. Its only pole is at zero, where
it has the same leading term as $\eta_0$. Consequently
\[
 \omega_{\chi_c}(r_h)-\eta_0=df
\]
for an odd polynomial $f$ in $t^{-1}$ of degree less than $K$.
As in the preceding proof, the single-factor assertion of
Lemma~\ref{lem-correction}, now at $r_h$, subtracts $df$. The
depths at $r_h$ and $r_c$ remain $Mr_h$ and $\alpha$, respectively,
so the slope between them remains $K$.

By construction, $\eta-\eta_0=dP_{\mathrm{aux}}$, where
$P_{\mathrm{aux}}$ is the sum of the principal parts of the
primitive at its nonzero poles. Equation~\eqref{eq-auxfactor}
gives
\[
 \overline{(G_h-1)/\varpi_{s_h}}=P_{\mathrm{aux}}
\]
at $r_h$. Hence $G_h$ adds precisely $\eta-\eta_0$ to the
differential by \eqref{eq-kummerdepth}. We have obtained a fibre
for which
\begin{equation}\label{eq-actual}
 \omega_\chi(r_h)=\eta,\qquad
 \delta_\chi(r_h)=Mr_h,\qquad
 \delta_\chi'(r_h^-)=M,
\end{equation}
where the last equality follows from
$\operatorname{ord}_\infty(\eta)=M-1$.

We next choose a fibre for which slope $M$ persists down to zero.
Let the roots of $G_c$ vary in $\mathcal A_c$. For $G_h$, vary
the principal-part coefficients in the closed unit disk, requiring
each leading coefficient to be a unit. If there is one pair of
centres, keep it fixed. If there are two pairs, keep the second
pair fixed and let the first centre $\xi$ vary subject to
\[
 |\xi|=|\xi^2-c_2|=|\xi+\xi_0|=1.
\]
Here $c_2$ lifts the square of the second centre and $\xi_0$ lifts
the chosen first centre. The middle condition keeps the two pairs
distinct. The last excludes the reduction $-\bar\xi_0$.

These conditions define an affinoid parameter space. By
Lemma~\ref{lem-auxiliary-branches}, a principal part of order
$j<p$ gives $j$ simple zeros and one pole, all reducing to its
centre. Thus the auxiliary factor supplies $2B$ branch points,
disjoint from the critical points since $r_h<r_c$. There are
$M+1$ branch points in all. The verification of connectedness,
pure inseparability and the closed disk condition in the previous
proof applies unchanged for $0<r\leq r_h$. We may therefore apply
Theorem~\ref{thm-minimum} with $r_0=r_h$ and slope bound $M$.

We need to know which fibres have left slope $M$ at $r_h$. We
show that every such fibre has differential $\eta$ and the
parameter reductions used in \eqref{eq-actual}. Write $L_s$ and
$R_s$ for the left and right slopes at $s$. The slope formulas
and the fact that a differential on $\mathbb P^1_k$ has divisor
of degree $-2$ give
\[
 L_s-R_s
 =-\sum_{z\ne0,\infty}\operatorname{ord}_z\omega_\chi(s).
\]
For $r_h<s<r_c$, no branch point reduces to any of the $z$ in
this sum. Hence $L_s\leq R_s$ by \cite[Proposition~6.19]{Obus}.
The slopes are therefore nondecreasing on this interval. They
are at most $K$, since exactly $K+1$ branch points reduce to zero.
At $s=r_h$, the nonzero finite directions contain $2B$ branch
points, so $L_{r_h}-R_{r_h}\leq2B$. If $L_{r_h}=M$, then
$R_{r_h}\geq M-2B=K$. It follows that the slope on $(r_h,r_c)$
equals $K$. Since $\delta_\chi(r_c)=\alpha$, we get
$\delta_\chi(r_h)=Mr_h$.

We first recover the critical root reductions. In
\eqref{eq-auxfactor}, the centres are units, so $\widehat A(u)$
has integral Taylor coefficients at $u=0$. Its constant term
cancels in $\widehat A(u)-\widehat A(-u)$. Taking
$u=T/\varpi_{r_h}$ gives
\[
 v_{r_c}(G_h-1)\geq s_h+r_c-r_h>0.
\]
Thus the character of $G_h$ has depth less than $\alpha$ at $r_c$
and does not change the differential there. The left slope $K$
at $r_c$ now gives the moment equations by \eqref{eq-criticalbase},
so the critical roots have the prescribed reductions. Moreover,
$\delta_{\chi_c}(r_h)\geq Mr_h$ by the slope bound $K$ for
$\chi_c$. Strict inequality is impossible, since the character
of $G_h$ has depth $Mr_h$ and could not cancel a term of larger
depth. Hence $\delta_{\chi_c}(r_h)=Mr_h$. The leading coefficient
of its differential at zero is therefore $C$ by
Lemma~\ref{lem-bounded-critical-depth}. The auxiliary differential
is regular at zero, so the full differential has this same
coefficient.

We now recover the auxiliary parameters. Since
$L_{r_h}-R_{r_h}=2B$, equality holds in each of the pole bounds
whose sum is $2B$. A centre containing $j+1$ branch points has a
pole of order $j+1$, and there are no other nonzero finite poles
or zeros. The orders at zero and infinity are $-K-1$ and $M-1$.
Thus the differential has the denominator in \eqref{eq-oneaux}
or \eqref{eq-twoaux} and has constant numerator. Once the centres
are known, the coefficient $C$ determines this numerator.

Only the first centre in the two-pair case remains to be
determined. Set $b'=\bar\xi^2$. Since $Mr_h>par_h$, the
differential is exact. Multiplying by $t^{pa}$ preserves exactness
because $d(t^{pa})=0$. On setting $z=t^2$, the resulting
differential is, up to a nonzero scalar, the pullback of
\[
 \frac{z^{a-1}\,dz}{(z-b')^{p-1}(z-c_2)^{b_2}}.
\]
This differential is exact over $k(z)$. Indeed, averaging a
primitive over $t\mapsto-t$ gives a primitive in $k(t^2)$.
Its residue at $c_2$ is therefore zero. Equation~\eqref{eq-residue}
gives $b'=ac_2/(a-b_2+1)=b$. The condition $|\xi+\xi_0|=1$ now
gives $\bar\xi=\bar\xi_0$. This proves that the differential is
$\eta$ in both cases. It also determines the principal-part
coefficients. Indeed, differentiation is injective on sums
$\sum_{j=1}^J a_j(t-\xi)^{-j}$ with $J<p$, since none of
$1,\ldots,J$ vanishes in $k$.

Choose a fibre minimizing
\[
 \mu=\sup\Bigl(\{0\}\cup
 \{r\in(0,r_h]\cap\mathbb Q\mid\delta_\chi'(r^-)<M\}\Bigr).
\]
The minimum exists by Theorem~\ref{thm-minimum}. The fibre in
\eqref{eq-actual} has slope $M$ on an interval ending at $r_h$,
so $\mu<r_h$. By what we have just proved, the minimizing fibre
has differential $\eta$ at $r_h$ and depth $Mr$ on $[\mu,r_h]$.

Suppose $\mu>0$. Then $\mu$ is a rational breakpoint. All branch
points have valuation greater than $\mu$, and $M\mu>pa\mu$.
Thus $\omega_\chi(\mu)$ is exact, and its only pole is at zero.
The right slope $M$ gives pole order $M+1$. Since $p\nmid M$,
we may write
\[
 \omega_\chi(\mu)=c_\mu t^{-M-1}dt+df
\]
with $c_\mu\ne0$ and $f$ an odd polynomial in $t^{-1}$ of
degree less than $M$. Apply Lemma~\ref{lem-correction} at $\mu$
with $\beta=\alpha-M\mu$ to a lift of $f$, taking the correction
to be $1$ if $df=0$. The first inequality in \eqref{eq-correction}
subtracts $df$ at $\mu$. The second preserves $\eta$ at $r_h$.
All parameter increments have positive valuation, so the corrected
fibre still belongs to the family.

The depths at $\mu$ and $r_h$ remain $M\mu$ and $Mr_h$. Since
the slopes are at most $M$, they equal $M$ between these endpoints.
The new differential at $\mu$ is $c_\mu t^{-M-1}dt$, so its left
slope is also $M$. The slope $M$ therefore extends below $\mu$,
contradicting minimality. Hence $\mu=0$. By continuity,
\begin{equation}\label{eq-goodline}
 \delta_\chi(r)=Mr\ \text{ for }0\leq r\leq r_h.
\end{equation}
In particular, $\delta_\chi(0)=0$ and $\delta_\chi'(0^+)=M$.
There are $a+1$ branch points of index $p^2$ and $N+2B=M-a$
of index $p$. Good reduction with jumps $(a,M)$ now follows
from \cite[Proposition~6.18]{Obus}.
\end{proof}

\begin{proof}[Proof of Proposition~\ref{prop-normalized}]
It remains to prescribe $f_2$. We first treat its leading
coefficient when $e>0$. In this case $M>pa$ and $p\nmid M$.
By \eqref{arb-jump-formula}, every cover in
Lemma~\ref{lem-bounded-both-depths} has a nonzero coefficient
$\gamma$ of $x^{-M}$ in its standard second coordinate.

Let $c$ denote the constant numerator of $\eta$ in
\eqref{eq-oneaux} or \eqref{eq-twoaux}. Multiply every squared
auxiliary centre by $\rho\in k^\times$, keeping $C$ fixed. The
relation between the two centres in \eqref{eq-residue} is
unchanged, while $c$ is replaced by $c\rho^B$. Thus every
$c\in k^\times$ occurs. We compute the corresponding change
in $\gamma$.

Take two such covers with numerators $c,c'$, coefficients
$\gamma,\gamma'$, and factors $G_1,G_2$. Write
\begin{equation}\label{eq-leaderraw}
 G_2/G_1=1+\sum_{i>0}b_iT^{-i}.
\end{equation}
We claim that $v(b_i)\geq\alpha$ for $i\geq M$, with strict
inequality for $i>M$. For a term of degree $i$, the bounds from
the two factors are
\[
 \begin{array}{c|c}
 \text{source of the term} & \text{valuation}\\ \hline
 \text{critical roots alone} & >ir_c\\
 \text{one auxiliary term} & \geq s_h+ir_h\\
 \text{at least two auxiliary terms} & \geq 2s_h+ir_h\\
 \text{an auxiliary term and critical degree }j>0
     & \geq s_h+ir_h+j(r_c-r_h).
 \end{array}
\]
The first bound follows from \eqref{eq-criticalraw}, since the
critical reductions are fixed. The remaining bounds follow by
expanding \eqref{eq-auxfactor}. The claim follows from
$s_h+Mr_h=\alpha<Mr_c$. These bounds also give convergence at
$r=0$.

In degree $M$, only the linear auxiliary term can contribute at
valuation $\alpha$. The coefficient of $t^{-M-1}dt$ in
$\eta'-\eta$ at infinity is $c'-c$. The principal parts at zero
have pole order at most $K+1<M+1$, so they make no contribution
in this degree. Hence
\[
 -M[t^{-M}](P'_{\mathrm{aux}}-P_{\mathrm{aux}})=c'-c.
\]
Returning from $t$ to $T$ multiplies the coefficient by
$\varpi_{s_h}\varpi_{r_h}^M=\varpi_\alpha$. With
$\lambda_p=\zeta_p-1$ and the normalization
$\varpi_\alpha=\lambda_p^p$, this gives
\[
 \overline{b_M/\lambda_p^p}=-(c'-c)/M.
\]

To identify this coefficient in the Witt vector, consider the
degree-$p$ character of $G_2/G_1$. It has depth zero since both
covers have depth zero at $r=0$. Apply Lemma~\ref{an-cutoff} with
$m=M$ and $\beta=\alpha$, taking the normalized representative
$F_1$ first and \eqref{eq-leaderraw} second. Then
$\overline{(F_1-1)/\lambda_p^p}$ is a polynomial in $x^{-1}$ of
degree at most $M$, with no terms of degree divisible by $p$.
Its coefficient in degree $M$ is $\overline{b_M/\lambda_p^p}$.
This polynomial is the right hand side of the residue
Artin--Schreier equation. Indeed, substituting $Y=1+\lambda_p Z$
in $Y^p=F_1$ and using $\overline{p/\lambda_p^{p-1}}=-1$ gives
\[
 Z^p-Z=\overline{(F_1-1)/\lambda_p^p}.
\]

The two Witt vectors have the same first coordinate. Their
difference is therefore $(0,f'_2-f_2)$, where $f_2$ and $f'_2$
here denote their standard second coordinates. Two standard
Artin--Schreier polynomials without constant term represent the
same class only if they are equal. Indeed, a nonzero pole of
$h^p-h$ has order divisible by $p$. Comparing coefficients gives
\begin{equation}\label{eq-leader}
 \gamma'-\gamma=-(c'-c)/M.
\end{equation}
Thus $\gamma+c/M$ has a value $A$ independent of the choice of
cover. If $A\ne0$, we could choose $c=MA\ne0$ and obtain
$\gamma=0$, contrary to the jump $M>pa$. Hence $\gamma=-c/M$.
In particular, we may prescribe any nonzero leading coefficient.

When $e=0$, no leading coefficient needs to be prescribed in
degree $M=pa$, since a standard coordinate has no term of degree
divisible by $p$. Lemma~\ref{lem-bounded-critical-depth} gives
$\delta_{\chi_c}(r)=Mr$ on $[0,r_c]$. There are $a+1$ branch
points of index $p^2$ and $(p-1)a=M-a$ of index $p$. Thus
\cite[Proposition~6.18]{Obus} gives good reduction with jumps
$(a,M)$ in this case as well.

In both cases, let $f$ be the prescribed second coordinate minus
the one just obtained. It is an odd standard polynomial in $x^{-1}$
of degree less than $M$. If $f\ne0$, take an odd lift $F$ of $-f$
and apply Lemma~\ref{lem-correction} at $r=0$ with $\beta=\alpha$.
Use only $G_c$ when $e=0$. The correction satisfies
\[
 \overline{(Q-1)/\lambda_p^p}=f.
\]
The substitution above shows that its residue equation is
$Z^p-Z=f$. Thus it adds the Witt vector $(0,f)$, giving precisely
the prescribed second coordinate. The branch counts and inertia
orders have not changed. Since the residue extension has jumps
$(a,M)$, the generic and special differents are equal by
\eqref{eq-budget}. The corrected cover therefore has good reduction.

Finally, let $\sigma$ generate the cyclic subgroup of order $p^2$
and let $\tau$ be the prescribed reflection. The equality of Witt
characters gives an isomorphism of residue extensions respecting
$\sigma$. Under this isomorphism the reduced reflection is
$\sigma^j\tau$ for some $j$. Composing the isomorphism with
$\sigma^i$ replaces $j$ by $j+2i$, since
$\tau\sigma\tau^{-1}=\sigma^{-1}$. As $p$ is odd, we may choose
$i$ with $j+2i\equiv0\pmod{p^2}$. The isomorphism then respects
both $\sigma$ and $\tau$, and hence identifies the reduction with
the given local $D_{p^2}$-action.
\end{proof}

\section{Proofs of the main theorems}\label{sec-proofs-main}
\label{sec-deform}\label{def-two-deformations}\label{def-composition}
\label{sec-arbitrary-length}\label{sec-characteristic-three}

We first prove Theorem~\ref{thm-main}. By
Proposition~\ref{prop-normalized}, a local $D_{p^2}$-action lifts
when its cyclic upper jumps are $(a,pa+2e)$ with $a<p$ and
$0\leq e<p$. To treat the remaining jumps, we deform the given
cover in characteristic $p$. On the generic fibre, the jumps at
the origin lie in the bounded range, and the other ramification
points have cyclic stabilizers. We use a similar deformation for
Theorem~\ref{thm-arbitrary-length}. The proof of
Theorem~\ref{thm-characteristic-three} instead uses an explicit
polynomial satisfying the lifting criterion of Obus.

Set $A=k[[\mu]]$, $K_\mu=k((\mu))$, and
$\Omega=\overline{K_\mu}$. Our deformations are given by Witt
equations over $K_\mu(x)$. Substituting $\mu=0$ will recover the
original equation. To show that we also recover the original
smooth curve, we compare the total degrees of the different.

\begin{lemma}\label{lem-genus}
Let $Y\to\mathbb P^1_k$ be a smooth connected cyclic cover of
degree $p^r$. Suppose that a Witt equation of length $r$ over
$K_\mu(x)$ is regular at the generic point of the special
$x$-line and reduces there to an equation for $k(Y)/k(x)$.
Then its geometric generic cover is connected. Let $\Delta_\eta$
and $\Delta_0$ be the total different degrees of this cover and of
$Y\to\mathbb P^1_k$, respectively. Then
$\Delta_\eta\geq\Delta_0$. If equality holds, the equation gives
a smooth family with special fibre $Y$ after a finite extension
of $A$.
\end{lemma}

\begin{proof}
Let $A'$ be the valuation ring of any finite extension of $K_\mu$.
Its residue field is still $k$. The local ring $B$ of
$\mathbb P^1_{A'}$ at the generic point of the special line is a
discrete valuation ring with residue field $k(x)$. The successive
Artin--Schreier equations are monic with derivative $-1$. They
therefore define a finite \'etale $B$-algebra of rank $p^r$.
Its reduction is the field $k(Y)$, so it is a local ring, and
hence a discrete valuation ring with the same uniformizer as $B$.
Its fraction field has degree $p^r$ over $\operatorname{Frac}(B)$.
This holds after every finite extension of the constants, proving
geometric connectedness.

Choose $A'$ so that the smooth projective model of the cover over
$\Omega$ and its map to the $x$-line are defined over
$\operatorname{Frac}(A')$. Let $\mathcal Y$ be the normalization
of $\mathbb P^1_{A'}$ in its function field, and let $C$ be the
special fibre. The normalization is finite since $A'$ is a
complete discrete valuation ring. It is flat over $A'$ since
it is torsion free.

Every component of $C$ maps onto the special line. Indeed, its
generic point has height one, and going down over the normal base
$\mathbb P^1_{A'}$ shows that its image has height one as well.
The calculation over $B$ therefore shows that $C$ has one
component, of multiplicity one, with function field $k(Y)$.
There are no embedded points, since $C$ is cut out by a
nonzerodivisor on a normal surface, whose local rings are
Cohen--Macaulay. Thus $C$ is integral and its normalization is $Y$.

Let $\nu\colon Y\to C$ be the normalization map and set
$Q=\nu_*\mathcal O_Y/\mathcal O_C$. Taking Euler characteristics
gives $p_a(C)-g(Y)=\dim_k H^0(C,Q)$, where $p_a(C)$ is the
arithmetic genus. See \cite[Tag~0CE0]{Stacks}. Properness and
flatness give $g(\mathcal Y_\Omega)=p_a(C)$. Riemann--Hurwitz
for the two covers of degree $p^r$ now gives
\[
 \Delta_\eta-\Delta_0
 =2\bigl(p_a(C)-g(Y)\bigr)
 =2\dim_k H^0(C,Q).
\]
If $\Delta_\eta=\Delta_0$, then $Q=0$ and $C=Y$. Both geometric
fibres are smooth, so the flat family $\mathcal Y$ is smooth
over $A'$.
\end{proof}

Thus it suffices to prove $\Delta_\eta\leq\Delta_0$. We obtain
this upper bound from the poles of the Witt coordinates. The
following lemma allows us to bound them before standardizing
the equation.

\begin{lemma}\label{lem-witt-pole-bounds}
Let $(f_1,\ldots,f_r)$ be a Witt vector over a Laurent series field
of characteristic $p$. Suppose that the pole order of $f_i$ is at
most $q_i$, where $q_1\geq0$ and $q_i\geq p q_{i-1}$ for
$2\leq i\leq r$. After a finite extension of the coefficient field,
the Witt class has a standard representative with the same bounds.
\end{lemma}

\begin{proof}
Give the variables in coordinate $j$ weight $p^{j-1}$. Every
monomial in coordinate $i$ of a Witt sum or difference has weight
$p^{i-1}$. If $m_j$ is its total exponent in the variables of
weight $p^{j-1}$, its pole order is at most
\[
 \sum_{j=1}^i m_jq_j
 \leq \frac{q_i}{p^{i-1}}\sum_{j=1}^i m_jp^{j-1}
 =q_i.
\]
Thus Witt addition and subtraction preserve the bounds.

Let $t$ be a local parameter. To remove a term $ct^{-pd}$ in
coordinate $j$, subtract $\wp(h)$, where $h_j=c^{1/p}t^{-d}$
and the other coordinates are zero. This replaces $ct^{-pd}$ by
$c^{1/p}t^{-d}$ and leaves the earlier coordinates unchanged.
Both $h$ and its Frobenius satisfy the given bounds, so the later
coordinates still satisfy them. Removing the largest divisible
pole order first makes the procedure finite. To remove a regular
part $g$, use a regular solution of $h_j^p-h_j=g$. After adjoining
a solution of the constant equation, Hensel's lemma gives $h_j$.
Proceeding from $j=1$ to $j=r$ proves the assertion.
\end{proof}

We shall also use the following consequence of the local--global
principle. It explains why it is enough to lift the local actions
on the geometric generic fibre of each deformation.

\begin{lemma}\label{lem-specialization}
Let $\mathcal Y\to\mathbb P^1_A$ be a smooth family of $G$-covers
with special fibre $Y\to\mathbb P^1_k$. If every local action at a
ramification point over $\Omega$ lifts over a finite extension of
$W(\Omega)$, then $Y\to\mathbb P^1_k$ lifts over a complete
discrete valuation ring finite over $W(k)$.
\end{lemma}

\begin{proof}
Let $y$ be a ramification point of the special fibre and let
$I=G_y$ be its stabilizer. Completing the family at $y$ and at
its image gives an extension $A[[s]]\subset A[[z]]$ whose
reduction is the given local $I$-action. The ramified local
actions on its geometric generic fibre come from the points of
$\mathcal Y_\Omega$ specializing to $y$, so they lift by
hypothesis. By \cite[Theorem~6.20]{ObusLifting}, the local action
at $y$ lifts over a finite extension of $W(k)$. Taking a common
finite extension for the finitely many ramified orbits, the
local--global principle \cite[Theorem~3.1]{ObusSurvey} gives the
required lift. We may take $\mathbb P^1_R$ as the lift of the
quotient over this common ring $R$.
\end{proof}

We now return to a local $D_{p^2}$-action with cyclic upper jumps
$(u,v)$. By Lemma~\ref{pre-normalized-witt}, its Witt vector may
be taken to be $(x^{-u},f_2)$, where $f_2$ is an odd standard
polynomial in $x^{-1}$. All different degrees below are computed
for the cyclic cover of the $x$-line. By \eqref{eq-different}, a
point of full inertia with jumps $(m_1,m_2)$ contributes
\[
 d(m_1,m_2)=(p-1)(m_1+1)+(p^2-p)(m_2+1).
\]
A branch point of inertia $p$ and jump $j$ contributes
$p(p-1)(j+1)$, since there are $p$ points above it. The special
cyclic cover is branched only at $0$, so $\Delta_0=d(u,v)$.

\begin{lemma}\label{lem-lastdeformation}
Suppose $u=a<p$ and $v\geq pa+2p$. Let $M$ be the unique integer
with $pa<M<pa+2p$ and $M\equiv v\pmod{2p}$, and set
$D=(v-M)/2$. There is a smooth deformation with cyclic upper
jumps $(a,M)$ at $0$. At each of $\mu$ and $-\mu$, the inertia
group has order $p$ and the jump is $D-1$.
\end{lemma}

\begin{proof}
Since $v>pa$, \eqref{arb-jump-formula} gives
$\deg_{x^{-1}}f_2=v$ and $p\nmid v$. Thus $M$ exists,
$p\nmid M$, and $D$ is a positive multiple of $p$. Set
$P(x)=x^vf_2(x)$. This is an even polynomial with $P(0)\ne0$
and $\deg P\leq v-1$. Keep the first coordinate fixed and take
\begin{equation}\label{eq-lastdeform}
 (F_1,F_2)=
 \left(x^{-a},\frac{P(x)}{x^M(x^2-\mu^2)^D}\right).
\end{equation}
Since $M+2D=v$, setting $\mu=0$ gives $(x^{-a},f_2)$. Both
coordinates are odd and regular away from $0,\mu,-\mu$, including
at infinity.

At $0$ the pole orders are $a$ and $M$. The first coordinate is
standard, and the leading term of the second survives
standardization since $p\nmid M$. As $M>pa$, the jumps are
$(a,M)$.

At either of $\mu$ and $-\mu$, the first coordinate is regular.
Over $\Omega$, choose a regular solution of $b^p-b=F_1$ in the
completed local ring. Subtracting $\wp(b,0)$ makes the first
coordinate zero and changes the second by a regular series.
The remaining Artin--Schreier equation has pole order at most
$D$. Since $p\mid D$, its jump, if it ramifies, is at most $D-1$.
Hence
\begin{equation}\label{def-last-budget}
 \Delta_\eta\leq d(a,M)+2p(p-1)D=d(a,v)=\Delta_0.
\end{equation}
Lemma~\ref{lem-genus} gives smoothness and equality. Each new
point must therefore contribute $p(p-1)D$, giving inertia $p$
and jump $D-1$ at both points.
\end{proof}

\begin{lemma}\label{lem-firstdeformation}
Suppose $u>p$. Set $D_1=(u-p)/2$ and $D_2=(v-p^2)/2$.
There is a smooth deformation with cyclic upper jumps
$(p-2,p^2-2)$ at $0$ and $(D_1,D_2)$ at each of $\mu$ and
$-\mu$. All three points have inertia $C_{p^2}$.
\end{lemma}

\begin{proof}
Since $u$ is odd and prime to $p$, the integer $D_1$ is positive
and prime to $p$. Also $D_2-pD_1=(v-pu)/2\geq0$. Set
$P(x)=x^vf_2(x)$, allowing $P=0$. Since
$\deg_{x^{-1}}f_2\leq v$, this is an even polynomial of degree
at most $v-1$ when nonzero. Take
\begin{equation}\label{eq-firstdeform2}
 F_1=\frac{1}{x^p(x^2-\mu^2)^{D_1}}
 \quad\text{and}\quad
 F_2=\frac{P(x)}{x^{p^2}(x^2-\mu^2)^{D_2}}.
\end{equation}
The identities $p+2D_1=u$ and $p^2+2D_2=v$ show that the
reduction is $(x^{-u},f_2)$. The only possible poles are
$0,\mu,-\mu$.

At each of $\mu$ and $-\mu$, the first coordinate has pole order
$D_1$. Thus the first quotient ramifies with jump $D_1$. The
inertia group is $C_{p^2}$, since no proper subgroup of $C_{p^2}$
maps onto its quotient $C_p$. Lemma~\ref{lem-witt-pole-bounds}
bounds the standard pole orders by $D_1,D_2$. Since $pD_1\leq D_2$,
the second jump is at most $D_2$.

At $0$ the pole orders are at most $p$ and $p^2$. We may
standardize within these bounds and preserve oddness, as in
Lemma~\ref{arb-standard-witt}. The largest odd integers prime
to $p$ in the two ranges are $p-2$ and $p^2-2$. If the first
quotient ramifies at $0$, then the inertia is $C_{p^2}$, and the
jumps are at most $p-2,p^2-2$, since $p(p-2)<p^2-2$.
If the first quotient does not ramify, its standard coordinate
is zero. Only an equation of degree at most $p$ remains. The
contribution above $0$ is then at most
\[
 p(p-1)(p^2-1)
 =d(p-2,p^2-2)-(p-1)^2
 <d(p-2,p^2-2).
\]

Adding the three bounds gives
\begin{equation}\label{def-first-budget}
 \Delta_\eta
 \leq d(p-2,p^2-2)+2d(D_1,D_2)
 =d(u,v)=\Delta_0.
\end{equation}
The middle equality follows from
$(p-1)+2(D_1+1)=u+1$ and $(p^2-1)+2(D_2+1)=v+1$.
Lemma~\ref{lem-genus} gives smoothness and equality. In
particular, the strict inequality above cannot hold, so the
inertia at $0$ is full. Both jumps enter the different with
positive coefficients. Hence each upper bound is attained at
each of the three points.
\end{proof}

\begin{proof}[Proof of Theorem~\ref{thm-main}]
If $u<p$ and $v<pu+2p$, then $v=pu+2e$ for some $0\leq e<p$,
and Proposition~\ref{prop-normalized} applies. Suppose otherwise.
If $u<p$, use Lemma~\ref{lem-lastdeformation}. If $u>p$, use
Lemma~\ref{lem-firstdeformation}. There is no case $u=p$, since
$p\nmid u$. Both deformations have odd Witt coordinates, so
$(x,y)\mapsto(-x,-y)$ extends to the normalization and reduces
to the prescribed reflection.

At $x=0$ the cyclic inertia is $C_{p^2}$, and the reflection
fixes the unique point above $0$. The local action there is
therefore dihedral. Its jumps are either $(u,M)$ with
$pu<M<pu+2p$, or $(p-2,p^2-2)$. Both pairs lie in the bounded
range, since $p^2-2=p(p-2)+2(p-1)$. The local action lifts by
Proposition~\ref{prop-normalized}.

The reflection exchanges $\mu$ and $-\mu$, so their stabilizers
are cyclic. At infinity the cyclic cover is unramified and
$s=x^2$ has index two. The cyclic Oort theorem
\cite[Theorem~1.1(1)]{Pop} gives lifts at these remaining points.
Lemma~\ref{lem-specialization} now gives a lift of the original
projective cover. Completing at the point above $s=0$ gives the
required action on $R[[Z]]$.
\end{proof}

We next consider a local $D_{p^n}$-action whose cyclic upper
jumps satisfy $u_i\equiv-1\pmod{2p}$ for every $i$. In this
case we can remove the cyclic ramification at $0$ entirely.
We use the different exponent of a totally ramified cyclic
extension of degree $p^n$,
\begin{equation}\label{arb-different}
 d(u_1,\ldots,u_n)=\sum_{i=1}^n(p^i-p^{i-1})(u_i+1).
\end{equation}
Indeed, set $u_0=-1$. On the $i$th interval the ramification
group has order $p^{n-i+1}$ and the corresponding length in
the lower numbering is $p^{i-1}(u_i-u_{i-1})$. Substituting
in Hilbert's different formula gives \eqref{arb-different}.

\begin{proposition}\label{arb-cyclic-satellite-deformation}
Suppose $u_i\equiv-1\pmod{2p}$ for $1\leq i\leq n$, and set
$D_i=(u_i+1)/2$. The projective $D_{p^n}$-cover admits a smooth
deformation in which the cyclic cover of the geometric generic
$x$-line is branched precisely at $\mu$ and $-\mu$. At both
points the inertia is $C_{p^n}$ and the jumps are
$D_1-1,\ldots,D_n-1$.
\end{proposition}

\begin{proof}
Choose a standard Witt representative $(f_1,\ldots,f_n)$ with
each $f_i$ odd. Since $p\nmid u_i$, the jump formula
\eqref{arb-jump-formula} gives $\deg_{x^{-1}}f_i=u_i$. Write
$P_i(x)=x^{u_i}f_i(x)$. Then $P_i$ is even, $P_i(0)\ne0$, and
$\deg P_i\leq u_i-1$. Take
\begin{equation}\label{arb-satellite-equations}
 F_i(x,\mu)=\frac{xP_i(x)}{(x^2-\mu^2)^{D_i}}
 \text{ for }1\leq i\leq n.
\end{equation}
Since $2D_i=u_i+1$, setting $\mu=0$ gives $f_i$. Each $F_i$ is
odd and is regular away from $\mu$ and $-\mu$, including at $0$
and infinity.

By Lemma~\ref{lem-genus}, the geometric generic cover is
connected. Let $I_+$ and $I_-$ denote the inertia subgroups at
$\mu$ and $-\mu$, allowing the trivial subgroup if the point is
unramified. The reflection exchanges these groups and acts on
$C_{p^n}$ by inversion, so $I_+=I_-$. Denote their common value
by $I$. The quotient by $I$ is a connected unramified cover of
$\mathbb P^1_\Omega$. If its degree is $m$, Riemann--Hurwitz
gives $2g-2=-2m$. Thus $m=1$ and $I=C_{p^n}$. By reflection
the two points have the same jumps, say $s_1,\ldots,s_n$.

We claim that $s_i\leq D_i-1$. For $i>1$, the difference
$u_i-pu_{i-1}$ is nonnegative and congruent to $p-1$ modulo $2p$.
Thus
\[
 D_i-pD_{i-1}
 =\frac{u_i-pu_{i-1}-(p-1)}{2}\geq0.
\]
The pole order of $F_i$ at $\mu$ is at most $D_i$.
Lemma~\ref{lem-witt-pole-bounds} preserves these bounds in
standard form. Since $p\mid D_i$, the standard pole orders are
at most $D_i-1$. In particular, $s_1\leq D_1-1$. If
$s_{i-1}\leq D_{i-1}-1$, then
$ps_{i-1}\leq pD_{i-1}-p<D_i-1$, so
\eqref{arb-jump-formula} gives $s_i\leq D_i-1$.

The special cyclic cover has one branch point and the generic
cyclic cover has two. Hence
\begin{equation}\label{arb-satellite-budget}
 \Delta_\eta
 =2\sum_{i=1}^n(p^i-p^{i-1})(s_i+1)
 \leq 2\sum_{i=1}^n(p^i-p^{i-1})D_i
 =\Delta_0.
\end{equation}
Lemma~\ref{lem-genus} gives smoothness and equality. Each
coefficient $p^i-p^{i-1}$ is positive, so $s_i=D_i-1$ for every
$i$. The reflection extends to the normalization and specializes
to the given reflection.
\end{proof}

\begin{proof}[Proof of Theorem~\ref{thm-arbitrary-length}]
On the geometric generic fibre, the points $x=\pm\mu$ lie above
$s=\mu^2$ and have inertia $C_{p^n}$. At $s=0$ and $s=\infty$
the cyclic cover is unramified and $s=x^2$ has index two, so the
inertia is $C_2$. These are all the branch points. Every local
action therefore lifts by the cyclic Oort theorem
\cite[Theorem~1.1(1)]{Pop}. Apply Lemma~\ref{lem-specialization}
and complete at the point above $s=0$.
\end{proof}

It remains to prove Theorem~\ref{thm-characteristic-three}.
Assume $p=3$. Set $a_0=0$ and define $a_i$ successively by
\[
 a_i\equiv u_i\pmod6
 \quad\text{and}\quad
 3a_{i-1}\leq a_i<3a_{i-1}+6.
\]
By \cite[Proposition~1.11]{Obus}, it suffices to lift all actions
with jumps $(a_1,\ldots,a_n)$ over every algebraically closed
field of characteristic three. When $a_i>3a_{i-1}$ and
$3\nmid a_{i-1}$ for $2\leq i\leq n$, the criterion of
\cite[Proposition~1.6 and Theorem~1.14]{Obus} takes the following
form. For each $U=a_{i-1}$, set $h=(U-1)/2$. We seek an even
polynomial $f(t)=\sum_{j=0}^{U-1}b_jt^{2j}$ of degree $2U-2$.
Write $\omega_f=dt/(t^{U+1}f(t))$ and set $b_j=0$ outside the
indicated range. By \cite[Definition~2.1 and Corollary~4.7]{DDKOT},
it is enough to verify
\begin{equation}\label{charthree-differential-criterion}
 \Cart(\omega_f)=\omega_f+Ut^{-U-1}dt
 \quad\text{and}\quad
 \det\bigl(b_{1+3i-j}\bigr)_{0\leq i,j<h}\ne0.
\end{equation}
When $h=0$, the determinant is $1$.

\begin{proposition}\label{charthree-explicit-data}
Let $U$ be positive with $U\equiv5\pmod6$, and set $h=(U-1)/2$.
Then
\begin{equation}\label{charthree-polynomial}
 f(t)=1+\sum_{j=h+1}^{2h}t^{2j}
\end{equation}
satisfies \eqref{charthree-differential-criterion}. In fact,
the matrix formed from the integer coefficients of $f$ has
determinant $1$ or $-1$.
\end{proposition}

\begin{proof}
Set $R(t)=t^U+t$. Since $U=-1$ in $k$, we have
\[
 (1-t^2)f(t)=1-t^2+t^{U+1}-t^{2U}=1-R(t)^2.
\]
Consequently,
\begin{equation}\label{charthree-logarithmic-identity}
 \omega_f-t^{-U-1}dt
 =\frac{t^{U-1}-1}{1-R^2}\,dt
 =\frac{2\,dR}{1-R^2}
 =d\log\left(\frac{1+R}{1-R}\right),
\end{equation}
where $d\log(g)=dg/g$. Cartier fixes $dg/g$ and kills
$t^{-U-1}dt=d(-t^{-U}/U)$. Applying Cartier gives
$\Cart(\omega_f)=\omega_f-t^{-U-1}dt$, as required. Also
$f(0)=1$. At a root of $f$, the left side of
\eqref{charthree-logarithmic-identity} has a pole whose order is
the multiplicity of the root. The right side has at most a simple
pole. Hence every root of $f$ is simple.

For the determinant, the only nonzero coefficients are $b_0=1$
and $b_j=1$ for $h+1\leq j\leq2h$. Work over $\mathbb F_\ell$
for an arbitrary prime $\ell$. Let $(c_0,\ldots,c_{h-1})$ be a
kernel vector and set $Y_j=\sum_{r=0}^j c_r$.

Index the rows by $s=1,4,\ldots,3h-2$. The entry in column $j$
is $b_{s-j}$, which equals $1$ precisely when $j=s$ or
$h+1\leq s-j\leq2h$. The row equations are therefore
\[
 0=\sum_{j=0}^{h-1}b_{s-j}c_j
 =
 \begin{cases}
 c_s & s<h,\\
 Y_{s-h-1} & h<s\leq2h,\\
 Y_{h-1}-Y_{s-2h-1} & s>2h.
 \end{cases}
\]
There is no row $s=h$, since $h\equiv2\pmod3$. The middle rows
give $Y_j=0$ for $j=1,4,\ldots,h-1$. The first rows give
$Y_j=Y_{j-1}$ at these same indices, so $Y_j=0$ also for
$j=0,3,\ldots,h-2$. In particular, $Y_{h-1}=0$. The last rows
now give $Y_j=0$ for the remaining indices $j=2,5,\ldots,h-3$.
This last list is empty when $h=2$. Thus all partial sums vanish
and every $c_j$ is zero. The integer determinant is nonzero modulo
every prime $\ell$, so it equals $1$ or $-1$.
\end{proof}

\begin{proof}[Proof of Theorem~\ref{thm-characteristic-three}]
For $n=1$ the assertion follows from Bouw and Wewers \cite{BW}.
Suppose $n\geq2$. The first jump is odd and prime to three, so
$a_1\in\{1,5\}$. For $2\leq i<n$, the hypothesis gives
$u_i\equiv5\pmod6$, whereas $3a_{i-1}\equiv3\pmod6$. Thus
$a_i=3a_{i-1}+2$. At the last step, $3\nmid u_n$ gives
$a_n-3a_{n-1}\in\{2,4\}$. Hence every step is strict, and each
preceding jump $U=a_{i-1}$ is either $1$ or congruent to $5$
modulo $6$.

If $U=1$, take $f=-1$. Then $\omega_f=-t^{-2}dt=d(t^{-1})$,
so both sides of the Cartier equation vanish and the determinant
is $1$. For $U\equiv5\pmod6$, take the polynomial of
Proposition~\ref{charthree-explicit-data}. These polynomials are
defined over $\mathbb F_3$, so the criterion holds over every
algebraically closed field of characteristic three. This proves
the theorem.
\end{proof}

\bibliographystyle{abbrv}
\bibliography{references}
\end{document}